\documentclass[9pt,a4paper]{amsart}%
\usepackage{amsfonts}
\usepackage[hiresbb]{graphicx}
\usepackage{amsmath}
\usepackage{amssymb}
\usepackage{version}
\usepackage{epsfig}
\usepackage{epstopdf}
\usepackage{bibmods}%
\providecommand{\U}[1]{\protect\rule{.1in}{.1in}}
\DeclareGraphicsExtensions{.eps}
\theoremstyle{theorem}
\newtheorem{defn}{Definition}[section]

\newtheorem{exmp}[defn]{Example}
\newtheorem{lem}[defn]{Lemma}

\newtheorem{prop}[defn]{Proposition}
\newtheorem{cor}[defn]{Corollary}

\newtheorem{notat}[defn]{Notation}

\newenvironment{pf}[1][Proof]{\textbf{#1.} }{\ \rule{0.5em}{0.5em}}

\begin{document}
\title[Bifurcation equations]{Bifurcation equations for periodic orbits of implicit discrete dynamical systems}
\author[H. Oliveira]{Henrique M. Oliveira}
\address{Av. rovisco Pais, 1049-001 Lisbon, Portugal\\
Center for Mathematical Analysis Geometry and Dynamical Systems, Mathematics
Department, Technical Institute of Lisbon, University of Lisbon}
\email{holiv@math.ist.ulisboa.pt}
\date{January, 2016 - AMS: Primary:37G15; Secondary: 39A28}

\begin{abstract}
Bifurcation equations, non-degeneracy and transversality conditions are
obtained for the fold, transcritical, pitchfork and flip bifurcations for
periodic points of one dimensional implicitly defined discrete dynamical
systems. The backward Euler method and the trapezoid method for numeric
solutions of ordinary differential equations fall in the category of implicit
dynamical systems. Examples of bifurcations are given for some implicit
dynamical systems including bifurcations for the backward Euler method when
the step size is changed.

\end{abstract}
\maketitle

\section{Introduction}

\subsection{Motivation}

In this paper we study bifurcation equations and transversality conditions for
local bifurcations of $p$-periodic points in one-dimensional discrete
dynamical systems defined implicitly, with $p$ a positive integer. In
particular, we focus our attention on the fold, transcritical, pitchfork and
flip, i.e., the most frequently found in applications. The main result of the
paper is to obtain expressions for the general bifurcation equations and
transversality conditions in dynamical systems defined implicitly.

Implicitly defined discrete and continuous dynamical systems are not very well
studied, only very recently Albert Luo published \textquotedblleft\textit{the
first monograph to discuss the implicit mapping dynamics of periodic flows to
chaos\textquotedblright} \cite{luo2015}. The singularities of some implicit
continuous dynamical systems in dimension two have been addressed in
\cite{Dav2008}, namely the Clairaut system. Nevertheless, it is an interesting
and open field of research. This type of dynamical system appears in
applications, namely in the theory of PDE in the works of Sharkovsky and
co-workers \cite{CE,LSS,SSSV1,SSSV2,SP}, in Mathematical Economics directly
\cite{Me2} or in the context of backward dynamics \cite{KSY,ME}. It appears
also in the context of Control Theory \cite{Holl2005}. These implicit
dynamical systems appear also in numerical methods for ordinary differential
equations, v.g., the backward Euler, the trapezoid method
\cite{Ush1986,Hir194} and the Runge-Kutta implicit method, see the recent
article \cite{Zhang2016}. Implicit numeric methods are very useful when the
original equations exhibit stiffness, see for instance \cite{Hai2009,Hai1996}.
In \cite{Lamarque2000} the implicit Euler method was used in a concrete
mechanical problem.  Some implicit iterative schemes were transformed  in
forward dynamical systems using numerical methods,  v.g., Newton method,
\cite{Ghe2003}. In implicit numerical schemes it is possible to prove the
existence of period doubling when the step size parameter increases as we do
with a simple example at the end of this article.  It is also interesting to
see the existence of chaos when the parameter $h$ is big enough, but still
relatively small.

The case of $p$-periodic points with $p>1$, is very intricate, the
computations increase its complexity extraordinary with the powers of the
normal form, as we can see in this paper in the case of the pitchfork. For
that reason, we study codimension $1$ cases, the most common in applications.

The study of one-dimensional bifurcations makes sense, since many higher order
systems can be reduced  \cite{Rega2005} to lower order dimensional dynamics
via center manifold and Poincar\'e map techniques as in \cite{Kle1995}, using
spectral  properties and quasi-periodicity \cite{Mezic2005}, and in periodic
non-autonomous systems using Floquet theory \cite{david2000}.

It is completely open and would be interesting to investigate the invariance
of the bifurcation equations for periodic non-autonomous systems defined
implicitly in the line of work of \cite{HO1}.

One of the main reasons of this paper is to provide computational tools for
the applied researcher dealing with implicitly defined dynamical systems. It
is possible to study the bifurcations that can occur without the knowledge of
an explicit difference equation. All the formulae are programmable using the
usual platforms available for mathematicians. The examples where prepared
using Wolfram Mathematica 10.0.

We follow the terminology of \cite{KU}.

\subsection{Overview}

We organized this paper in four sections. In Section \ref{SectionPrelim} we
introduce basic concepts.

In Section 3, the core of this work, we study in detail the equations of
bifurcation for $p$-periodic orbits of implicitly defined dynamical systems.

In Section 4 we present examples, namely on the Euler method for numerical
solutions of ordinary differential equations. In the implicit difference
equations of numerical methods we show the existence of bifurcation depending
on the step size parameter $h$, and the existence of chaos even in very simple examples.

\section{\label{SectionPrelim}Preliminaries}

\subsection{Basic definitions and notation}

We define implicitly a discrete dynamical system using instead of the classic
definition%
\begin{equation}
x_{n+1}=f(x_{n})\text{, }x_{n}\in I\text{, with }n\in\mathbb{%
%TCIMACRO{\U{2115} }%
%BeginExpansion
\mathbb{N}
%EndExpansion
}\text{,} \label{pp}%
\end{equation}
the alternative one%
\[
F\left(  x_{n},x_{n+1}\right)  =0\text{, for }x_{n},x_{+1}\in I\text{, with
}n\in\mathbb{%
%TCIMACRO{\U{2115} }%
%BeginExpansion
\mathbb{N}
%EndExpansion
}\text{,}%
\]
where $I$ is a real interval (not necessarily compact and maybe $%
%TCIMACRO{\U{211d} }%
%BeginExpansion
\mathbb{R}
%EndExpansion
$), where we input $x_{n}$ and solve for $x_{n+1}$ giving an initial condition
$x_{0}\in I$. The usual Euclidean distance is defined in $I$. The map $F$ is
sufficiently differentiable for the purposes of bifurcation theory, assumption
that we keep in this paper. We suppose that given $F\left(  x,y\right)  =0$,
there exists the solution $\left(  x_{0},y_{0}\right)  $, and an implicit
function $y=f\left(  x\right)  $ with $y_{0}=f\left(  x_{0}\right)  $ such
that%
\[
F\left(  x,f\left(  x\right)  \right)  =0\text{, }%
\]
in a suitable neighborhood of $\left(  x_{0},y_{0}\right)  $. We follow
\cite{KP} concerning the implicit function theorem. For the purposes of this
article we admit the existence of the necessary solutions in the appropriate
neighborhoods of the bifurcation points. Obviously, each particular dynamical
system defined implicitly must be studied to ensure the existence of the
iteration function $f\left(  x\right)  $.

In the sequel, by ${\mathcal{C}}\left(  D\right)  $ we denote the collection
of all continuous maps in its domain $D$, by ${\mathcal{C}}^{1}\left(
D\right)  $ the collection of all continuously differentiable elements of
${\mathcal{C}}\left(  D\right)  $ and, in general by ${\mathcal{C}}^{s}\left(
D\right)  ,$ $s\geq1,$ the collection of all elements of ${\mathcal{C}}\left(
D\right)  $ having continuous derivatives up to order $s$ in $D.$

The $p$ composition of $f$ a real function of real variable is denoted by
$f^{p}$, the usual power is denoted by $\left(  f\right)  ^{p}$.

Let $f\in{\mathcal{C}}^{1}\left(  D\right)  $, and let $x_{0}$ be a periodic
point of period $p$, $x_{0}$ is called a \emph{hyperbolic attractor} if
$|\frac{df^{p}(x_{0})}{dx}|<1$, a \emph{hyperbolic repeller} if $|\frac
{df^{p}(x_{0})}{dx}|>1$, and \emph{non-hyperbolic} if $|\frac{df^{p}(x_{0}%
)}{dx}|=1$.

\begin{defn}
We say that two continuous maps $f:I\rightarrow I$ and $g:J\rightarrow J$, are
topologically conjugate, if there exists a homeomorphism $h:I\rightarrow J$,
such that $h\circ f=g\circ h$. We call $h$ the topological conjugacy of $f$
and $g$.
\end{defn}

We use $\alpha$ for a real parameter.

\begin{defn}
\label{BifSet}If $f\left(  \cdot,\alpha\right)  $ is a family of maps, then
the regular values $\alpha$ of the parameters are those which have the
property that $f\left(  \cdot,\widetilde{\alpha}\right)  $ is topologically
conjugate to $f\left(  \cdot,\alpha\right)  $ for all $\widetilde{\alpha}$ in
some open neighbourhood of $\alpha$. If $\alpha$ is not a regular value, it is
a \emph{bifurcation value}. The collection of all the bifurcation values is
the \emph{bifurcation set}, $\Omega\subset\mathbb{R}$, in the parameter space.
\end{defn}

Let $f\left(  \cdot,\alpha_{0}\right)  $ be a parameter dependent family of
maps in ${\mathcal{C}}^{s}\left(  D\right)  $. Let $\alpha_{0}$ be a
particular parameter and $a\in D$ be a fixed point of the $p$ composition map
$f\left(  \cdot,\alpha_{0}\right)  $, with $p$ a minimal positive integer,
i.e.,%
\[
a=f^{p}\left(  a,\alpha_{0}\right)  \text{,}%
\]
$a$ is a periodic point of the dynamical system. The condition of $a$ being
non-hyperbolic is necessary for the existence of a local bifurcation. The
existence and nature of that bifurcation depends on other symmetry and
differentiable conditions that we will see bellow. If there exists a local
bifurcation we say that $(a,\alpha_{0})$ is a \emph{bifurcation point }(when
there is no risk of confusion, we say that $a$ is a \emph{bifurcation point}).

\begin{notat}
For notational simplicity we consider the real parameter $\alpha$ as a
standard variable along with the dynamic variable $x$, i.e., we write
$F\left(  x,y,\alpha\right)  $ instead of $F_{\alpha}\left(  x,y\right)  $,
reserving the last slot for the parameter, keeping in mind that the
compositions are always in the dynamic variables $x$ and $y$. In this paper we
never use $f_{\alpha}$ to mean dependence on the parameter.

When there is no danger of confusion and no operations regarding the
parameter, we denote the evaluation of functions depending on the dynamic
variable and the parameter omitting the later, for instance $F\left(
x,y,\alpha\right)  $ or $f\left(  x,\alpha\right)  $ will be denoted by
$F\left(  x,y\right)  $ or $f\left(  x\right)  $ in order to avoid to overload
the complicated notation needed for the computations of chain rules.
Nevertheless, all the maps in this paper depend on the parameter as well on
the dynamic variable. We deal with parameter depending families of maps, even
when that dependence is not visible in some formulas or expressions.

We denote the derivatives relative to some variable $y$ by $\partial_{y}$.
Repeated differentiation relative to the same variable is denoted by
$\partial_{y^{n}}$, for instance $\partial_{yyy}=\partial_{y^{3}}$. When there
is no danger of confusion, we denote strict partial derivatives, i.e., not
seeing composed functions, by a subscript. For instance, the third partial
derivative of $f$ relative to $y$ is, in that case, denoted by $f_{yyy}$ or
$f_{y^{3}}$.
\end{notat}

This means, in particular, that when dealing with the composition of real
scalar functions $F\left(  x,y\right)  $ with $g\left(  x,y\right)  $ and
$h\left(  x,y\right)  $, we have the usual chain rule
\begin{multline*}
\partial_{x}F\left(  g\left(  x,t\right)  ,h\left(  x,t\right)  ,\alpha
\right)  =\\
F_{x}\left(  g\left(  x,t\right)  ,h\left(  x,t\right)  ,\alpha\right)
g_{x}\left(  x,t\right)  +F_{y}\left(  g\left(  x,t\right)  ,h\left(
x,t\right)  ,\alpha\right)  h_{x}\left(  x,t\right)  \text{,}%
\end{multline*}

\subsection{Classic conditions for fold, transcritical, pitchfork and flip
bifurcations}

In this paragraph, we recall briefly the conditions of codimension $1$ local
bifurcations with derivatives $\partial_{x}f^{p}(x_{0})=\pm1$.

We first consider the case $\partial_{x}f^{p}(x_{0})=+1$. Giving a discrete
dynamical system generated by the iteration of $f$ in its domain $D$, and a
real parameter $\alpha$, in order to compute the bifurcation points one has to
solve the bifurcation equations \cite{KU}%
\begin{equation}%
\begin{array}
[c]{l}%
{f}^{p}(x,\alpha)=x\text{, fixed point equation}\\
{f}_{x}^{p}(x,\alpha)=1\text{, non-hyperbolicity condition.}%
\end{array}
\label{Fold}%
\end{equation}

\subsubsection{Fold}

The simplest of such local bifurcations is the fold or saddle node
bifurcation. One assumes, in this case, the non-degeneracy condition
\begin{equation}
{f_{x^{2}}}(x,\alpha)\not =0 \label{Nondeg}%
\end{equation}
and the transversality condition \cite{KU}%
\begin{equation}
f_{\alpha}\left(  x,\alpha\right)  \not =0. \label{unfold1}%
\end{equation}
We set generically that $\alpha\in%
%TCIMACRO{\U{211d} }%
%BeginExpansion
\mathbb{R}
%EndExpansion
$, since one needs only one parameter to unfold locally this singularity
\cite{A,AR,CH,GS,GU,KU}. The normalized germ of this bifurcation is%
\[
x\pm x^{2}\text{,}%
\]
with principal family%
\[
x\pm x^{2}+\alpha\text{,}%
\]
which is locally weak topologically conjugated to any other family
\cite{AR,KU} satisfying the bifurcation conditions.

\subsubsection{Transcritical}

Another simple bifurcation is the transcritical, in this case is a bifurcation
with symmetry. One assumes, in this case, the non-degeneracy condition
\begin{equation}
{f_{x^{2}}}(x,\alpha)\not =0,
\end{equation}
the transversality condition of the fold fails%
\begin{equation}
f_{\alpha}\left(  x,\alpha\right)  =0,
\end{equation}
becoming a new degeneracy condition. The symmetry condition states that the
fixed point of $f$ persists. Without loss of generality we consider that $0$
is that fixed point. The new transversality condition is%
\[
f_{x\alpha}\left(  x,\alpha\right)  \not =0\text{.}%
\]
Again, we set generically that $\alpha\in%
%TCIMACRO{\U{211d} }%
%BeginExpansion
\mathbb{R}
%EndExpansion
$, since one needs only one parameter to unfold locally this singularity
\cite{A,AR,CH,GS,GU}. The principal family is now%
\[
\left(  1+\alpha\right)  x\pm x^{2}\text{,}%
\]
which is weak topologically conjugated to any other family \cite{AR,KU}
satisfying the bifurcation conditions.

\subsubsection{Pitchfork}

The last type of bifurcation we consider with derivative $\partial_{x}%
f^{p}(x_{0})=+1$ is the pitchfork, another bifurcation with the same symmetry
on the fixed point as the transcritical. One assumes, in this case, the extra
degeneracy condition%
\[
{f_{x^{2}}}(x,\alpha)=0,
\]
and the new non-degeneracy condition
\begin{equation}
{f_{x^{3}}}(x,\alpha)\not =0.
\end{equation}
The transversality condition of the fold fails again%
\begin{equation}
f_{\alpha}\left(  x,\alpha\right)  =0\text{,}%
\end{equation}
and the transversality condition is assumed again to be%
\[
f_{x\alpha}\left(  x,\alpha\right)  \not =0\text{.}%
\]
We set generically that $\alpha\in%
%TCIMACRO{\U{211d} }%
%BeginExpansion
\mathbb{R}
%EndExpansion
$ \cite{A,AR,CH,GS,GU}. The principal family is now%
\[
\left(  1+\alpha\right)  x\pm x^{3}\text{,}%
\]
which is weak topologically conjugated to any other family \cite{AR,KU}
satisfying the bifurcation conditions.

\subsubsection{Flip}

We consider now the conditions of codimension $1$ local bifurcations with
derivative $\partial_{x}f^{p}(x_{0})=-1$.

One has to solve the bifurcation equations \cite{KU}%
\begin{equation}%
\begin{array}
[c]{l}%
{f}^{p}(x,\alpha)=x\text{, fixed point equation}\\
{f}_{x}^{p}(x,\alpha)=-1\text{, non-hyperbolicity condition.}%
\end{array}
\label{Flip}%
\end{equation}
One assumes, in this case, the generic non-degeneracy condition
\begin{equation}
\frac{1}{2}\left(  f_{x^{2}}\left(  x,\alpha\right)  \right)  ^{2}+\frac{1}%
{3}f_{x^{3}}\left(  x,\alpha\right)  \not =0, \label{Schwarz}%
\end{equation}
which is equivalent to say that the Schwarzian derivative
\[
Sf\left(  x,\alpha\right)  =\frac{f_{x^{3}}\left(  x,\alpha\right)  }%
{f_{x}\left(  x,\alpha\right)  }-\frac{3}{2}\left(  \frac{f_{x^{2}}\left(
x,\alpha\right)  }{f_{x}\left(  x,\alpha\right)  }\right)  ^{2}%
\]
of $f$ is not zero at the bifurcation point where $f_{x}\left(  x,\alpha
\right)  =-1$. The transversality condition \cite{KU} is%
\begin{equation}
f_{x\alpha}(0,0)\not =0\text{.} \label{fliptransverse}%
\end{equation}
We set generically that $\alpha\in%
%TCIMACRO{\U{211d} }%
%BeginExpansion
\mathbb{R}
%EndExpansion
$ \cite{A,AR,CH,GS,GU,KU}. The normalized germ of this bifurcation is%
\[
-x\pm x^{3}\text{,}%
\]
with principal family%
\[
-\left(  1+\alpha\right)  x\pm x^{3}\text{,}%
\]
which is again locally weak topologically conjugated to any other family
\cite{AR,KU} satisfying the bifurcation conditions.

Adding degeneracy conditions, one obtains higher degeneracy (higher
codimension) local bifurcations. In this paper we keep it simple and do not
consider higher codimension.

\section{Implicit discrete dynamical systems}

\subsection{Bifurcation equations}

Let us now consider the case of implicit DDS. Given the parameter depend
family $F\in{\mathcal{C}}^{s}\left(
%TCIMACRO{\U{211d} }%
%BeginExpansion
\mathbb{R}
%EndExpansion
^{2}\right)  $, with $s\geq1$, enough for our results, such that
\[%
\begin{array}
[c]{cccc}%
F: &
%TCIMACRO{\U{211d} }%
%BeginExpansion
\mathbb{R}
%EndExpansion
^{2} & \longrightarrow &
%TCIMACRO{\U{211d} }%
%BeginExpansion
\mathbb{R}
%EndExpansion
\text{,}\\
& \left(  x,y\right)  & \longmapsto & F\left(  x,y\right)  \text{.}%
\end{array}
\]
We start by the example of dynamics near fixed points. So, consider $F\left(
x_{f},x_{f}\right)  =0$, with derivative $F_{y}\left(  x_{f},x_{f}\right)
\not =0$. We have the implicit discrete dynamical system near the fixed point
$\left(  x_{f},x_{f}\right)  $ defined by%
\[
F\left(  x_{n},x_{n+1},\alpha\right)  =0\text{, for }x_{n},x_{+1}\in I\text{,
with }n\in\mathbb{%
%TCIMACRO{\U{2115} }%
%BeginExpansion
\mathbb{N}
%EndExpansion
}\text{.}%
\]
Along this work we always consider the independent variable in the first slot
of $F\left(  \cdot,\cdot,\cdot\right)  $, being the dependent variable, or
implicit function, at the second slot and the parameter at the third slot. One
instance of this type of systems is obtained by Sharkovsky and coauthors
\cite{CE,LSS,SSSV1,SSSV2,SP} in some boundary value problems. The classic
counterpart of this scheme is%
\[
x_{n+1}-f\left(  x_{n},\alpha\right)  =0\text{, for }x_{n},x_{+1}\in I\text{,
with }n\in\mathbb{%
%TCIMACRO{\U{2115} }%
%BeginExpansion
\mathbb{N}
%EndExpansion
}\text{,}%
\]
with a fixed point $x_{f}$ and with $F\left(  x,y,\alpha\right)  =y-f\left(
x,\alpha\right)  $. The classic bifurcation equations are relative to
$y=f\left(  x,\alpha\right)  $. The bifurcation equations in the implicit case
are%
\begin{align*}
F\left(  x,y\left(  x\right)  ,\alpha\right)   &  =0\text{,}\\
F_{x}\left(  x,y\left(  x\right)  ,\alpha\right)  +F_{y}\left(  x,y\left(
x\right)  ,\alpha\right)  y_{x}\left(  x\right)   &  =0\text{.}%
\end{align*}
At the bifurcation point $y=x=x_{f}$, we have $y_{x}\left(  x_{f}\right)
=f_{x}\left(  x_{f},\alpha\right)  =\pm1$, the equations become%
\begin{align*}
F\left(  x_{f},x_{f},\alpha\right)   &  =0\text{,}\\
F_{x}\left(  x_{f},x_{f},\alpha\right)  \pm F_{y}\left(  x_{f},x_{f}%
,\alpha\right)   &  =0\text{,}%
\end{align*}
with non-degeneracy condition%
\[
f_{x^{2}}\left(  x_{f},\alpha\right)  =-\frac{F_{x^{2}}\left(  x_{f}%
,x_{f},\alpha\right)  \pm2F_{xy}\left(  x_{f},x_{f},\alpha\right)  +F_{y^{2}%
}\left(  x_{f},x_{f},\alpha\right)  }{F_{y}\left(  x_{f},x_{f},\alpha\right)
}\not =0\text{.}%
\]

The case of periodic points is more involved, the orbit of $x$ is obtained by
successive substitution at the function $F\left(  x,y,\alpha\right)  $,
accordingly to the scheme%
\begin{equation}
\left\{
\begin{array}
[c]{l}%
F\left(  x,f\left(  x\right)  \right)  =0,\\
F\left(  f\left(  x\right)  ,f^{2}\left(  x\right)  \right)  =0,\\
\cdots\\
F\left(  f^{j-2}\left(  x\right)  ,f^{j-1}\left(  x\right)  \right)  =0,\\
F\left(  f^{j-1}\left(  x\right)  ,f^{j}\left(  x\right)  \right)  =0,\\
\cdots
\end{array}
\right.  \label{iterative}%
\end{equation}
or, with initial condition $x_{0}$%
\begin{equation}
\left\{
\begin{array}
[c]{l}%
F\left(  x_{0},x_{1}\right)  =0,\\
F\left(  x_{1},x_{2}\right)  =0,\\
\cdots\\
F\left(  x_{j-2},x_{j-1}\right)  =0,\\
F\left(  x_{j-1},x_{j}\right)  =0,\\
\cdots
\end{array}
\right.
\end{equation}
where we omitted $\alpha$ for the sake notational simplicity. In this case, we
suppose that there exists an implicit solution of $F\left(  x,y\right)  =0$,
such that $y=f\left(  x\right)  $ is well defined for all the points $x_{0},$
$\ldots,$ $x_{j}$, $\ldots$ meaning that $F_{y}\left(  x_{0},x_{1}\right)
\not =0$, $F_{y}\left(  x_{1},x_{2}\right)  \not =0$, $\ldots$, $F_{y}\left(
x_{p-1},x_{0}\right)  \not =0$, $\ldots$.

Naturally, $x_{0}$ is a periodic point of the implicit dynamical system if%
\begin{equation}
\left\{
\begin{array}
[c]{l}%
F\left(  x_{0},x_{1}\right)  =0,\\
F\left(  x_{1},x_{2}\right)  =0,\\
\cdots\\
F\left(  x_{p-2},x_{p-1}\right)  =0,\\
F\left(  x_{p-1},x_{0}\right)  =0\text{.}%
\end{array}
\right.  \label{implic}%
\end{equation}

To obtain the bifurcation equations for periodic points we compute the
derivatives of the system (\ref{implic}). The next two lemmas \ref{chain rule}
and \ref{chainparameter} are fundamental in the study of the bifurcation
conditions, giving explicit formulas for the computation of derivatives
relative to $x$ and the parameter $\alpha$. All the other derivatives used in
this paper in the bifurcation conditions whatsoever are obtained recursively
using the results of this two lemmas. The next Lemma \ref{chain rule}%
\ establishes the chain rule for the first derivative of $f$, the iteration
function defined implicitly by $F\left(  x,y\right)  =0$.

\begin{lem}
\label{chain rule}Chain rule for implicit orbits. The derivative of $f^{j}$
defined using the system (\ref{implic}) is given by
\begin{equation}
\partial_{x}f^{j}\left(  x\right)  =\left(  -1\right)  ^{j}%
%TCIMACRO{\dprod \limits_{i=0}^{j-1}}%
%BeginExpansion
{\displaystyle\prod\limits_{i=0}^{j-1}}
%EndExpansion
\frac{F_{x}\left(  f^{i}\left(  x\right)  ,f^{i+1}\left(  x\right)  \right)
}{F_{y}\left(  f^{i}\left(  x\right)  ,f^{i+1}\left(  x\right)  \right)  }.
\label{chain}%
\end{equation}
Equivalently, given the initial condition $x_{0}$
\begin{equation}
\partial_{x}f^{j}\left(  x_{0}\right)  =\left(  -1\right)  ^{j}%
%TCIMACRO{\dprod \limits_{i=0}^{j-1}}%
%BeginExpansion
{\displaystyle\prod\limits_{i=0}^{j-1}}
%EndExpansion
\frac{F_{x}\left(  x_{i},x_{i+1}\right)  }{F_{y}\left(  x_{i},x_{i+1}\right)
}. \label{chain2}%
\end{equation}

\end{lem}

\begin{pf}
We differentiate the system (\ref{iterative}) relative to $x$, noticing that
the zeroth order composition is the identity $f^{0}\left(  x\right)  =x$, and
$f_{x}^{0}\left(  x\right)  =1$, with the simplifying notation $f^{j}\left(
x\right)  =f^{j}$ for $j=0,1,2\ldots.$%
\begin{align*}
F_{x}\left(  f^{0},f\right)  +F_{y}\left(  f^{0},f\right)  f_{x}\left(
f^{0}\right)   &  =0,\\
F_{x}\left(  f,f^{2}\right)  f_{x}\left(  f^{0}\right)  +F_{y}\left(
f,f^{2}\right)  f_{x}\left(  f^{0}\right)  f_{x}\left(  f\right)   &  =0,\\
&  \cdots\\
F_{x}\left(  f^{j-1},f^{j}\right)
%TCIMACRO{\dprod \limits_{i=0}^{j-2}}%
%BeginExpansion
{\displaystyle\prod\limits_{i=0}^{j-2}}
%EndExpansion
f_{x}\left(  f^{i}\right)  +F_{y}\left(  f^{j-1},f^{j}\right)
%TCIMACRO{\dprod \limits_{i=0}^{j-1}}%
%BeginExpansion
{\displaystyle\prod\limits_{i=0}^{j-1}}
%EndExpansion
f_{x}\left(  f^{i}\right)   &  =0,\\
&  \cdots
\end{align*}
cancelling the common factors we get%
\begin{align*}
F_{x}\left(  f^{0},f\right)  +F_{y}\left(  f^{0},f\right)  f_{x}\left(
f^{0}\right)   &  =0,\\
F_{x}\left(  f,f^{2}\right)  +F_{y}\left(  f,f^{2}\right)  f_{x}\left(
f\right)   &  =0,\\
&  \cdots\\
F_{x}\left(  f^{j-1},f^{j}\right)  +F_{y}\left(  f^{j-1},f^{j}\right)
f_{x}\left(  f^{j-1}\right)   &  =0,\\
&  \cdots
\end{align*}
solving for $f_{x}\left(  f^{j}\right)  $ we obtain%
\begin{align*}
f_{x}\left(  f^{0}\right)   &  =-\frac{F_{x}\left(  f^{0},f\right)  }%
{F_{y}\left(  f^{0},f\right)  },\\
f_{x}\left(  f\right)   &  =-\frac{F_{x}\left(  f,f^{2}\right)  }{F_{y}\left(
f,f^{2}\right)  },\\
&  \cdots\\
f_{x}\left(  f^{j-1}\right)   &  =-\frac{F_{x}\left(  f^{j-1},f^{j}\right)
}{F_{y}\left(  f^{j-1},f^{j}\right)  },\\
&  \cdots.
\end{align*}
Using the chain rule along the orbit, one obtains the product
\[
\partial_{x}f^{j}\left(  x\right)  =%
%TCIMACRO{\dprod \limits_{i=0}^{j-1}}%
%BeginExpansion
{\displaystyle\prod\limits_{i=0}^{j-1}}
%EndExpansion
f_{x}\left(  f^{i}\right)  =\left(  -1\right)  ^{j}%
%TCIMACRO{\dprod \limits_{i=0}^{j-1}}%
%BeginExpansion
{\displaystyle\prod\limits_{i=0}^{j-1}}
%EndExpansion
\frac{F_{x}\left(  f^{i},f^{i+1}\right)  }{F_{y}\left(  f^{i},f^{i+1}\right)
}.
\]
The second relation (\ref{chain2}) is a simple reformulation of the first one
(\ref{chain}).
\end{pf}

\begin{cor}
We have the first bifurcation equation%
\begin{subequations}
\begin{align}
\partial_{x}f^{p}\left(  x_{0}\right)   &  =%
%TCIMACRO{\dprod \limits_{j=0}^{p-1}}%
%BeginExpansion
{\displaystyle\prod\limits_{j=0}^{p-1}}
%EndExpansion
f_{x}\left(  x_{j}\right)  =\pm1\label{a}\\
&  =\left(  -1\right)  ^{p}%
%TCIMACRO{\dprod \limits_{j=0}^{p-1}}%
%BeginExpansion
{\displaystyle\prod\limits_{j=0}^{p-1}}
%EndExpansion
\frac{F_{x}\left(  x_{j},x_{j+1\left(  \operatorname{mod}p\right)  }\right)
}{F_{y}\left(  x_{j},x_{j+1\left(  \operatorname{mod}p\right)  }\right)  }%
=\pm1. \label{b}%
\end{align}

\end{subequations}
\end{cor}

\begin{pf}
We consider that $x_{p}=x_{0}$ and substitute in the chain rule (\ref{chain2})
of Lemma \ref{chain rule}. The non-hyperbolicity condition is $\partial
_{x}f^{p}\left(  x_{0}\right)  =\pm1$.
\end{pf}

To decide if there is a bifurcation and its type is necessary to obtain the
transversality conditions using the parameter derivative. The first possible
condition involves $\partial_{\alpha}f^{p}$. The next Lemma
\ref{chainparameter} is fundamental in that concern.

\begin{lem}
\label{chainparameter}The derivative of $f^{j}$ relative to the parameter
$\alpha$ defined using the system (\ref{implic}) is given by
\begin{equation}
\partial_{\alpha}f^{j}=\left(  -1\right)  ^{j}%
%TCIMACRO{\dsum \limits_{k=0}^{j-1}}%
%BeginExpansion
{\displaystyle\sum\limits_{k=0}^{j-1}}
%EndExpansion
\frac{\left(  -1\right)  ^{k}F_{\alpha}\left(  f^{k},f^{k+1}\right)  }%
{F_{y}\left(  f^{k},f^{k+1}\right)  }%
%TCIMACRO{\dprod \limits_{i>k}^{j-1}}%
%BeginExpansion
{\displaystyle\prod\limits_{i>k}^{j-1}}
%EndExpansion
\frac{F_{x}\left(  f^{i},f^{i+1}\right)  }{F_{y}\left(  f^{i},f^{i+1}\right)
}. \label{derivalfa}%
\end{equation}

\end{lem}

\begin{pf}
Similar to the proof of Lemma \ref{chain rule}. We have now the general rule%
\[
\partial_{\alpha}F\left(  g,h\right)  =F_{x}\left(  g,h\right)  \partial
_{\alpha}g+F_{y}\left(  g,h\right)  \partial_{\alpha}h+F_{\alpha}\left(
g,h\right)  =0,
\]
the first derivative is
\[
F_{y}\left(  f^{0},f\right)  f_{\alpha}+F_{\alpha}\left(  f^{0},f\right)  =0,
\]
solving for $f_{\alpha}$%
\[
f_{\alpha}=-\frac{F_{\alpha}\left(  f^{0},f\right)  }{F_{y}\left(
f^{0},f\right)  }.
\]
Doing the same for the second composition we obtain%
\[
\partial_{\alpha}f^{2}=\frac{F_{\alpha}\left(  f^{0},f\right)  F_{x}\left(
f,f^{2}\right)  }{F_{y}\left(  f^{0},f\right)  F_{y}\left(  f,f^{2}\right)
}-\frac{F_{\alpha}\left(  f,f^{2}\right)  }{F_{y}\left(  f,f^{2}\right)  },
\]
for the third composition%
\[
\partial_{\alpha}f^{3}=-\frac{F_{\alpha}\left(  f^{0},f\right)  F_{x}\left(
f,f^{2}\right)  F_{x}\left(  f^{2},f^{3}\right)  }{F_{y}\left(  f^{0}%
,f\right)  F_{y}\left(  f,f^{2}\right)  F_{y}\left(  f^{2},f^{3}\right)
}+\frac{F_{\alpha}\left(  f,f^{2}\right)  F_{x}\left(  f^{2},f^{3}\right)
}{F_{y}\left(  f,f^{2}\right)  F_{y}\left(  f^{2},f^{3}\right)  }%
-\frac{F_{\alpha}\left(  f^{2},f^{3}\right)  }{F_{y}\left(  f^{2}%
,f^{3}\right)  }.
\]
The previous expressions suggest the general formula for the derivatives
relative to $\alpha$%
\[
\partial_{\alpha}f^{k}=\left(  -1\right)  ^{k}%
%TCIMACRO{\dsum \limits_{j=0}^{k-1}}%
%BeginExpansion
{\displaystyle\sum\limits_{j=0}^{k-1}}
%EndExpansion
\frac{\left(  -1\right)  ^{j}F_{\alpha}\left(  f^{j},f^{j+1}\right)  }%
{F_{y}\left(  f^{j},f^{j+1}\right)  }%
%TCIMACRO{\dprod \limits_{i>j}^{k-1}}%
%BeginExpansion
{\displaystyle\prod\limits_{i>j}^{k-1}}
%EndExpansion
\frac{F_{x}\left(  f^{i},f^{i+1}\right)  }{F_{y}\left(  f^{i},f^{i+1}\right)
},
\]
which is the induction hypothesis. Consider the general formula
\[
\partial_{\alpha}F\left(  f^{k},f^{k+1}\right)  =F_{x}\left(  f^{k}%
,f^{k+1}\right)  \partial_{\alpha}f^{k}+F_{y}\left(  f^{k},f^{k+1}\right)
\partial_{\alpha}f^{k+1}+F_{\alpha}\left(  f^{k},f^{k+1}\right)  =0,
\]
solving for $\partial_{\alpha}f^{k+1}$ we have%
\[
\partial_{\alpha}f^{k+1}=\frac{-F_{x}\left(  f^{k},f^{k+1}\right)
\partial_{\alpha}f^{k}-F_{\alpha}\left(  f^{k},f^{k+1}\right)  }{F_{y}\left(
f^{k},f^{k+1}\right)  }%
\]%
\[
=\scriptstyle\frac{-F_{x}\left(  f^{k},f^{k+1}\right)  \left(  \left(
-1\right)  ^{k}%
%TCIMACRO{\dsum \limits_{j=0}^{k-1}}%
%BeginExpansion
{\displaystyle\sum\limits_{j=0}^{k-1}}
%EndExpansion
\frac{\left(  -1\right)  ^{j}F_{\alpha}\left(  f^{j},f^{j+1}\right)  }%
{F_{y}\left(  f^{j},f^{j+1}\right)  }%
%TCIMACRO{\dprod \limits_{i>j}^{k-1}}%
%BeginExpansion
{\displaystyle\prod\limits_{i>j}^{k-1}}
%EndExpansion
\frac{F_{x}\left(  f^{i},f^{i+1}\right)  }{F_{y}\left(  f^{i},f^{i+1}\right)
}\right)  }{F_{y}\left(  f^{k},f^{k+1}\right)  }-\frac{F_{\alpha}\left(
f^{k},f^{k+1}\right)  }{F_{y}\left(  f^{k},f^{k+1}\right)  }%
\]%
\[
=\scriptstyle\left(  -1\right)  ^{k+1}%
%TCIMACRO{\dsum \limits_{j=0}^{k-1}}%
%BeginExpansion
{\displaystyle\sum\limits_{j=0}^{k-1}}
%EndExpansion
\frac{\left(  -1\right)  ^{j}F_{\alpha}\left(  f^{j},f^{j+1}\right)  }%
{F_{y}\left(  f^{j},f^{j+1}\right)  }%
%TCIMACRO{\dprod \limits_{i>j}^{k-1}}%
%BeginExpansion
{\displaystyle\prod\limits_{i>j}^{k-1}}
%EndExpansion
\frac{F_{x}\left(  f^{i},f^{i+1}\right)  F_{x}\left(  f^{k},f^{k+1}\right)
}{F_{y}\left(  f^{i},f^{i+1}\right)  F_{y}\left(  f^{k},f^{k+1}\right)
}+\left(  -1\right)  ^{2k+1}\frac{F_{\alpha}\left(  f^{k},f^{k+1}\right)
}{F_{y}\left(  f^{k},f^{k+1}\right)  }%
\]%
\[
=\left(  -1\right)  ^{k+1}%
%TCIMACRO{\dsum \limits_{j=0}^{k}}%
%BeginExpansion
{\displaystyle\sum\limits_{j=0}^{k}}
%EndExpansion
\frac{\left(  -1\right)  ^{j}F_{\alpha}\left(  f^{j},f^{j+1}\right)  }%
{F_{y}\left(  f^{j},f^{j+1}\right)  }%
%TCIMACRO{\dprod \limits_{i>j}^{k}}%
%BeginExpansion
{\displaystyle\prod\limits_{i>j}^{k}}
%EndExpansion
\frac{F_{x}\left(  f^{i},f^{i+1}\right)  }{F_{y}\left(  f^{i},f^{i+1}\right)
},
\]
as desired.
\end{pf}

\begin{cor}
In particular, at the bifurcation point the derivative relative to the
parameter takes the form%
\begin{equation}
\partial_{\alpha}f^{p}\left(  x_{0}\right)  =\left(  -1\right)  ^{p}%
%TCIMACRO{\dsum \limits_{j=0}^{p-1}}%
%BeginExpansion
{\displaystyle\sum\limits_{j=0}^{p-1}}
%EndExpansion
\frac{\left(  -1\right)  ^{j}F_{\alpha}\left(  x_{j},x_{j+1}\right)  }%
{F_{y}\left(  x_{j},x_{j+1}\right)  }%
%TCIMACRO{\dprod \limits_{i>j}^{p-1}}%
%BeginExpansion
{\displaystyle\prod\limits_{i>j}^{p-1}}
%EndExpansion
\frac{F_{x}\left(  x_{i},x_{i+1}\right)  }{F_{y}\left(  x_{i},x_{i+1}\right)
}\text{.} \label{Dalphabifpoint}%
\end{equation}

\end{cor}

To obtain the non-degeneracy conditions we have to compute the second
derivative of $F$. In the next proposition we obtain an explicit expression
for the second derivative.

For the next results we introduce the notation $F\left(  f^{j},f^{j+1}\right)
=F^{j}$, $F_{x}\left(  f^{j},f^{j+1}\right)  =F_{x}^{j}$, $F_{x^{2}}\left(
f^{j},f^{j+1}\right)  =F_{x^{2}}^{j}$, $F_{y}\left(  f^{j},f^{j+1}\right)
=F_{y}^{j}$, $F_{y^{2}}\left(  f^{j},f^{j+1}\right)  =F_{y^{2}}^{j}$,
$F_{xy}\left(  f^{j},f^{j+1}\right)  =F_{xy}^{j}$, $F_{\alpha}\left(
f^{j},f^{j+1}\right)  =F_{\alpha}^{j}$, $F_{x\alpha}\left(  f^{j}%
,f^{j+1}\right)  =F_{x\alpha}^{j}$, $F_{y\alpha}\left(  f^{j},f^{j+1}\right)
=F_{y\alpha}^{j}$. At the bifurcation point we use the notation $F\left(
x_{j},x_{j+1}\right)  =\widetilde{F}^{j}$, $F_{x}\left(  x_{j},x_{j+1}\right)
=\widetilde{F}_{x}^{j}$, $F_{x^{2}}\left(  x_{j},x_{j+1}\right)
=\widetilde{F}_{x^{2}}^{j}$, $F_{y}\left(  x_{j},x_{j+1}\right)
=\widetilde{F}_{y}^{j}$, $F_{y^{2}}\left(  x_{j},x_{j+1}\right)
=\widetilde{F}_{y^{2}}^{j}$, $F_{xy}\left(  x_{j},x_{j+1}\right)
=\widetilde{F}_{xy}^{j}$, $F_{\alpha}\left(  x_{j},x_{j+1}\right)
=\widetilde{F}_{\alpha}^{j}$, $F_{x\alpha}\left(  x_{j},x_{j+1}\right)
=\widetilde{F}_{x\alpha}^{j}$, $F_{y\alpha}\left(  x_{j},x_{j+1}\right)
=\widetilde{F}_{y\alpha}^{j}$ and the abbreviation
\[
\nu_{j}=\frac{F_{x}^{j}}{F_{y}^{j}}\text{ and }\widetilde{\nu}_{j}%
=\frac{\widetilde{F}_{x}^{j}}{\widetilde{F}_{y}^{j}}.
\]

\begin{prop}
\label{secd}The second derivative of $f^{k}$ defined using the system
(\ref{implic}) along the orbit is%
\begin{equation}
\partial_{x^{2}}f^{k}=\partial_{x}f^{k}%
%TCIMACRO{\dsum \limits_{j=0}^{k-1}}%
%BeginExpansion
{\displaystyle\sum\limits_{j=0}^{k-1}}
%EndExpansion
\frac{F_{x^{2}}^{j}-2F_{xy}^{j}\nu_{j}+F_{y^{2}}^{j}\nu_{j}^{2}}{F_{x}^{j}%
}\partial_{x}f^{j}. \label{SEC}%
\end{equation}
At the bifurcation point where $\partial_{x}f^{p}\left(  x_{0}\right)  =\pm1$,
with $x_{p}=x_{0}$, the second derivative takes the form%
\begin{equation}
\partial_{x^{2}}f^{p}\left(  x_{0}\right)  =\pm%
%TCIMACRO{\dsum \limits_{j=0}^{p-1}}%
%BeginExpansion
{\displaystyle\sum\limits_{j=0}^{p-1}}
%EndExpansion
\frac{\widetilde{F}_{x^{2}}^{j}-2\widetilde{F}_{xy}^{j}\widetilde{\nu}%
_{j}+\widetilde{F}_{y^{2}}^{j}\widetilde{\nu}_{j}^{2}}{\widetilde{F}_{x}^{j}%
}\partial_{x}f^{j}.
\end{equation}

\end{prop}

\begin{pf}
We recall (\ref{chain})%
\[
\partial_{x}f^{k}=\left(  -1\right)  ^{k}%
%TCIMACRO{\dprod \limits_{j=0}^{k-1}}%
%BeginExpansion
{\displaystyle\prod\limits_{j=0}^{k-1}}
%EndExpansion
\frac{F_{x}^{j}}{F_{y}^{j}}.
\]

The second derivative is%
\[
\partial_{x^{2}}f^{k}=%
%TCIMACRO{\dsum \limits_{j=0}^{k-1}}%
%BeginExpansion
{\displaystyle\sum\limits_{j=0}^{k-1}}
%EndExpansion
\left(  -1\right)  ^{k}\partial_{x}F_{x}^{j}%
%TCIMACRO{\dprod \limits_{i\not =j=0}^{k-1}}%
%BeginExpansion
{\displaystyle\prod\limits_{i\not =j=0}^{k-1}}
%EndExpansion
\frac{F_{x}^{i}}{F_{y}^{i}}-%
%TCIMACRO{\dsum \limits_{j=0}^{k-1}}%
%BeginExpansion
{\displaystyle\sum\limits_{j=0}^{k-1}}
%EndExpansion
\left(  -1\right)  ^{k}\frac{\partial_{x}F_{y}^{j}}{\left(  F_{y}^{j}\right)
^{2}}%
%TCIMACRO{\dprod \limits_{i\not =j=0}^{k-1}}%
%BeginExpansion
{\displaystyle\prod\limits_{i\not =j=0}^{k-1}}
%EndExpansion
\frac{F_{x}^{i}}{F_{y}^{i}},
\]
i.e.,%
\[
\partial_{x^{2}}f^{k}=\scriptstyle%
%TCIMACRO{\dsum \limits_{j=0}^{k-1}}%
%BeginExpansion
{\displaystyle\sum\limits_{j=0}^{k-1}}
%EndExpansion
\frac{F_{x^{2}}^{j}\partial_{x}f^{j}+F_{xy}^{j}\partial_{x}f^{j+1}}{F_{x}^{j}%
}\left(  -1\right)  ^{k}%
%TCIMACRO{\dprod \limits_{i=0}^{k-1}}%
%BeginExpansion
{\displaystyle\prod\limits_{i=0}^{k-1}}
%EndExpansion
\frac{F_{x}^{i}}{F_{y}^{i}}-%
%TCIMACRO{\dsum \limits_{j=0}^{k-1}}%
%BeginExpansion
{\displaystyle\sum\limits_{j=0}^{k-1}}
%EndExpansion
\frac{F_{xy}^{j}\partial_{x}f^{j}+F_{y^{2}}^{j}\partial_{x}f^{j+1}}{F_{y}^{j}%
}\left(  -1\right)  ^{k}%
%TCIMACRO{\dprod \limits_{i=0}^{k-1}}%
%BeginExpansion
{\displaystyle\prod\limits_{i=0}^{k-1}}
%EndExpansion
\frac{F_{x}^{i}}{F_{y}^{i}},
\]
which is%
\begin{align*}
\partial_{x^{2}}f^{k}  &  =%
%TCIMACRO{\dsum \limits_{j=0}^{k-1}}%
%BeginExpansion
{\displaystyle\sum\limits_{j=0}^{k-1}}
%EndExpansion
\left(  \frac{F_{x^{2}}^{j}\partial_{x}f^{j}+F_{xy}^{j}\partial_{x}f^{j+1}%
}{F_{x}^{j}}-\frac{F_{xy}^{j}\partial_{x}f^{j}+F_{y^{2}}^{j}\partial
_{x}f^{j+1}}{F_{y}^{j}}\right)  \partial_{x}f^{k}\\
&  =%
%TCIMACRO{\dsum \limits_{j=0}^{k-1}}%
%BeginExpansion
{\displaystyle\sum\limits_{j=0}^{k-1}}
%EndExpansion
\left(  \frac{F_{x^{2}}^{j}+F_{xy}^{j}\partial_{x}f^{j+1}}{F_{x}^{j}}%
-\frac{F_{xy}^{j}\partial_{x}f^{j}+F_{y^{2}}^{j}\partial_{x}f^{j+1}}{F_{y}%
^{j}}\right)  \partial_{x}f^{k},
\end{align*}
substituting in the above expression the values of $\partial_{x}f^{j}$ and
$\partial_{x}f^{j+1}$, such that
\[
\partial_{x}f^{j}=\left(  -1\right)  ^{j}%
%TCIMACRO{\dprod \limits_{i=0}^{j-1}}%
%BeginExpansion
{\displaystyle\prod\limits_{i=0}^{j-1}}
%EndExpansion
\frac{F_{x}^{i}}{F_{y}^{i}}=\left(  -1\right)  ^{j}%
%TCIMACRO{\dprod \limits_{i=0}^{j-1}}%
%BeginExpansion
{\displaystyle\prod\limits_{i=0}^{j-1}}
%EndExpansion
\nu_{i}%
\]
and
\[
\partial_{x}f^{j+1}=\left(  -1\right)  ^{j+1}%
%TCIMACRO{\dprod \limits_{i=0}^{j}}%
%BeginExpansion
{\displaystyle\prod\limits_{i=0}^{j}}
%EndExpansion
\frac{F_{x}^{i}}{F_{y}^{i}}=\left(  -1\right)  ^{j}%
%TCIMACRO{\dprod \limits_{i=0}^{j-1}}%
%BeginExpansion
{\displaystyle\prod\limits_{i=0}^{j-1}}
%EndExpansion
\nu_{i},
\]
we obtain
\begin{align*}
\partial_{x^{2}}f^{k}  &  =\partial_{x}f^{k}%
%TCIMACRO{\dsum \limits_{j=0}^{k-1}}%
%BeginExpansion
{\displaystyle\sum\limits_{j=0}^{k-1}}
%EndExpansion
\frac{\left(  -1\right)  ^{j}}{F_{x}^{j}}\left(
%TCIMACRO{\dprod \limits_{i=0}^{j-1}}%
%BeginExpansion
{\displaystyle\prod\limits_{i=0}^{j-1}}
%EndExpansion
\frac{F_{x}^{i}}{F_{y}^{i}}\right)  \left(  F_{x^{2}}^{j}-2F_{xy}^{j}%
\frac{F_{x}^{j}}{F_{y}^{j}}+F_{y^{2}}^{j}\left(  \frac{F_{x}^{j}}{F_{y}^{j}%
}\right)  ^{2}\right) \\
&  =\partial_{x}f^{k}%
%TCIMACRO{\dsum \limits_{j=0}^{k-1}}%
%BeginExpansion
{\displaystyle\sum\limits_{j=0}^{k-1}}
%EndExpansion
\frac{F_{x^{2}}^{j}-2F_{xy}^{j}\nu_{j}+F_{y^{2}}^{j}\nu_{j}^{2}}{F_{x}^{j}%
}\partial_{x}f^{j},
\end{align*}
as desired.

The second statement is immediate.
\end{pf}

The mixed derivative $\partial_{\alpha x}f^{p}$ is also necessary for some
computations in the case of transcritical, pitchfork and flip.

\begin{prop}
At the bifurcation point we have%
\[
\partial_{\alpha x}f^{p}\left(  x_{0}\right)  =\scriptstyle\pm%
%TCIMACRO{\dsum \limits_{j=0}^{p-1}}%
%BeginExpansion
{\displaystyle\sum\limits_{j=0}^{p-1}}
%EndExpansion
\left(  \frac{\widetilde{F}_{x^{2}}^{j}-2\widetilde{F}_{xy}^{j}\widetilde{\nu
}_{j}+\widetilde{F}_{y^{2}}^{j}\widetilde{\nu}_{j}^{2}}{\widetilde{F}_{x}^{j}%
}\partial_{\alpha}f^{j}+\left(  \frac{\widetilde{F}_{y^{2}}^{j}}{\widetilde
{F}_{y}^{j}}-\frac{\widetilde{F}_{xy}^{j}}{\widetilde{F}_{x}^{j}}\right)
\frac{\widetilde{F}_{\alpha}^{j}}{\widetilde{F}_{y}^{j}}+\frac{\widetilde
{F}_{x\alpha}^{j}}{\widetilde{F}_{x}^{j}}-\frac{\widetilde{F}_{y\alpha}^{j}%
}{\widetilde{F}_{y}^{j}}\right)  ,
\]
where
\[
\partial_{\alpha}f^{j}=\left(  -1\right)  ^{j}%
%TCIMACRO{\dsum \limits_{k=0}^{j-1}}%
%BeginExpansion
{\displaystyle\sum\limits_{k=0}^{j-1}}
%EndExpansion
\frac{\left(  -1\right)  ^{k}\widetilde{F}_{\alpha}^{k}}{\widetilde{F}_{y}%
^{k}}%
%TCIMACRO{\dprod \limits_{i>j}^{p-1}}%
%BeginExpansion
{\displaystyle\prod\limits_{i>j}^{p-1}}
%EndExpansion
\widetilde{\nu}_{i}.
\]

\end{prop}

\begin{pf}
We have now the derivative of (\ref{chain}),%
\begin{align*}
\partial_{\alpha x}f^{p} &  =\left(  -1\right)  ^{p}\partial_{\alpha}\left(
\frac{%
%TCIMACRO{\dprod \limits_{j=0}^{p-1}}%
%BeginExpansion
{\displaystyle\prod\limits_{j=0}^{p-1}}
%EndExpansion
F_{x}^{j}}{%
%TCIMACRO{\dprod \limits_{j=0}^{p-1}}%
%BeginExpansion
{\displaystyle\prod\limits_{j=0}^{p-1}}
%EndExpansion
F_{y}^{j}}\right)  \\
&  =\left(  -1\right)  ^{p}\frac{\partial_{\alpha}\left(
%TCIMACRO{\dprod \limits_{j=0}^{p-1}}%
%BeginExpansion
{\displaystyle\prod\limits_{j=0}^{p-1}}
%EndExpansion
F_{x}^{j}\right)  }{%
%TCIMACRO{\dprod \limits_{j=0}^{p-1}}%
%BeginExpansion
{\displaystyle\prod\limits_{j=0}^{p-1}}
%EndExpansion
F_{y}^{j}}-\left(  -1\right)  ^{p}\frac{\partial_{\alpha}\left(
%TCIMACRO{\dprod \limits_{j=0}^{p-1}}%
%BeginExpansion
{\displaystyle\prod\limits_{j=0}^{p-1}}
%EndExpansion
F_{y}^{j}\right)
%TCIMACRO{\dprod \limits_{j=0}^{p-1}}%
%BeginExpansion
{\displaystyle\prod\limits_{j=0}^{p-1}}
%EndExpansion
F_{x}^{j}}{\left(
%TCIMACRO{\dprod \limits_{j=0}^{p-1}}%
%BeginExpansion
{\displaystyle\prod\limits_{j=0}^{p-1}}
%EndExpansion
F_{y}^{j}\right)  ^{2}},
\end{align*}
at the bifurcation point we have%
\[
\left(  -1\right)  ^{p}\frac{%
%TCIMACRO{\dprod \limits_{j=0}^{p-1}}%
%BeginExpansion
{\displaystyle\prod\limits_{j=0}^{p-1}}
%EndExpansion
\widetilde{F}_{x}^{j}}{%
%TCIMACRO{\dprod \limits_{j=0}^{p-1}}%
%BeginExpansion
{\displaystyle\prod\limits_{j=0}^{p-1}}
%EndExpansion
\widetilde{F}_{y}^{j}}=\pm1.
\]
Therefore,%
\[
\partial_{\alpha x}f^{p}=\frac{\left(  -1\right)  ^{p}\partial_{\alpha}\left(
%
%TCIMACRO{\dprod \limits_{j=0}^{p-1}}%
%BeginExpansion
{\displaystyle\prod\limits_{j=0}^{p-1}}
%EndExpansion
\widetilde{F}_{x}^{j}\right)  \mp\partial_{\alpha}\left(
%TCIMACRO{\dprod \limits_{j=0}^{p-1}}%
%BeginExpansion
{\displaystyle\prod\limits_{j=0}^{p-1}}
%EndExpansion
\widetilde{F}_{y}^{j}\right)  }{%
%TCIMACRO{\dprod \limits_{j=0}^{p-1}}%
%BeginExpansion
{\displaystyle\prod\limits_{j=0}^{p-1}}
%EndExpansion
\widetilde{F}_{y}^{j}},
\]%
\[
=\frac{\left(  -1\right)  ^{p}%
%TCIMACRO{\dsum \limits_{j=0}^{p-1}}%
%BeginExpansion
{\displaystyle\sum\limits_{j=0}^{p-1}}
%EndExpansion
\partial_{\alpha}\widetilde{F}_{x}^{j}\left(
%TCIMACRO{\dprod \limits_{i\not =j=0}^{p-1}}%
%BeginExpansion
{\displaystyle\prod\limits_{i\not =j=0}^{p-1}}
%EndExpansion
\widetilde{F}_{x}^{i}\right)  \mp%
%TCIMACRO{\dsum \limits_{j=0}^{p-1}}%
%BeginExpansion
{\displaystyle\sum\limits_{j=0}^{p-1}}
%EndExpansion
\partial_{\alpha}\widetilde{F}_{y}^{j}\left(
%TCIMACRO{\dprod \limits_{i\not =j=0}^{p-1}}%
%BeginExpansion
{\displaystyle\prod\limits_{i\not =j=0}^{p-1}}
%EndExpansion
\widetilde{F}_{y}^{i}\right)  }{%
%TCIMACRO{\dprod \limits_{j=0}^{p-1}}%
%BeginExpansion
{\displaystyle\prod\limits_{j=0}^{p-1}}
%EndExpansion
\widetilde{F}_{y}^{j}}%
\]
and at the bifurcation point this is%
\begin{align*}
\partial_{\alpha x}f^{p} &  =\pm%
%TCIMACRO{\dsum \limits_{j=0}^{p-1}}%
%BeginExpansion
{\displaystyle\sum\limits_{j=0}^{p-1}}
%EndExpansion
\frac{\partial_{\alpha}\widetilde{F}_{x}^{j}}{\widetilde{F}_{x}^{j}}\mp%
%TCIMACRO{\dsum \limits_{j=0}^{p-1}}%
%BeginExpansion
{\displaystyle\sum\limits_{j=0}^{p-1}}
%EndExpansion
\frac{\partial_{\alpha}\widetilde{F}_{y}^{j}}{\widetilde{F}_{y}^{j}}\\
&  =\pm%
%TCIMACRO{\dsum \limits_{j=0}^{p-1}}%
%BeginExpansion
{\displaystyle\sum\limits_{j=0}^{p-1}}
%EndExpansion
\left(  \frac{\widetilde{F}_{x^{2}}^{j}\partial_{\alpha}f^{j}+\widetilde
{F}_{xy}^{j}\partial_{\alpha}f^{j+1}+\widetilde{F}_{x\alpha}^{j}}%
{\widetilde{F}_{x}^{j}}-\frac{\widetilde{F}_{xy}^{j}\partial_{\alpha}%
f^{j}+\widetilde{F}_{y^{2}}^{j}\partial_{\alpha}f^{j+1}+\widetilde{F}%
_{y\alpha}^{j}}{\widetilde{F}_{y}^{j}}\right)  .
\end{align*}
Knowing that%
\begin{align*}
\partial_{\alpha}f^{j+1} &  =\left(  -1\right)  ^{j+1}%
%TCIMACRO{\dsum \limits_{k=0}^{j}}%
%BeginExpansion
{\displaystyle\sum\limits_{k=0}^{j}}
%EndExpansion
\frac{\left(  -1\right)  ^{k}F_{\alpha}^{k}}{F_{y}^{k}}%
%TCIMACRO{\dprod \limits_{i>k}^{j}}%
%BeginExpansion
{\displaystyle\prod\limits_{i>k}^{j}}
%EndExpansion
\frac{F_{x}^{i}}{F_{y}^{i}}\\
&  =-\frac{F_{x}^{j}\partial_{\alpha}f^{j}+F_{\alpha}^{j}}{F_{y}^{j}},
\end{align*}
we have%
\begin{align*}
\partial_{\alpha x}f^{p} &  =\\
&  =\scriptstyle\pm%
%TCIMACRO{\dsum \limits_{j=0}^{p-1}}%
%BeginExpansion
{\displaystyle\sum\limits_{j=0}^{p-1}}
%EndExpansion
\left(  \frac{\widetilde{F}_{x^{2}}^{j}\partial_{\alpha}f^{j}-\widetilde
{F}_{xy}^{j}\frac{\widetilde{F}_{x}^{j}\partial_{\alpha}f^{j}+\widetilde
{F}_{\alpha}^{j}}{\widetilde{F}_{y}^{j}}+\widetilde{F}_{x\alpha}^{j}%
}{\widetilde{F}_{x}^{j}}-\frac{\widetilde{F}_{xy}^{j}\partial_{\alpha}%
f^{j}-\widetilde{F}_{y^{2}}^{j}\frac{\widetilde{F}_{x}^{j}\partial_{\alpha
}f^{j}+\widetilde{F}_{\alpha}^{j}}{\widetilde{F}_{y}^{j}}+\widetilde
{F}_{y\alpha}^{j}}{\widetilde{F}_{y}^{j}}\right)  \\
&  =\scriptstyle\pm%
%TCIMACRO{\dsum \limits_{j=0}^{p-1}}%
%BeginExpansion
{\displaystyle\sum\limits_{j=0}^{p-1}}
%EndExpansion
\left(  \frac{\left(  \widetilde{F}_{x^{2}}^{j}\frac{\widetilde{F}_{y}^{j}%
}{\widetilde{F}_{x}^{j}}+\widetilde{F}_{y^{2}}^{j}\frac{\widetilde{F}_{x}^{j}%
}{\widetilde{F}_{y}^{j}}-2\widetilde{F}_{xy}^{j}\right)  \partial_{\alpha
}f^{j}}{\widetilde{F}_{y}^{j}}+\left(  \frac{\widetilde{F}_{y^{2}}^{j}%
}{\widetilde{F}_{y}^{j}}-\frac{\widetilde{F}_{xy}^{j}}{\widetilde{F}_{x}^{j}%
}\right)  \widetilde{F}_{\alpha}^{j}+\left(  \frac{\widetilde{F}_{x\alpha}%
^{j}}{\widetilde{F}_{x}^{j}}-\frac{\widetilde{F}_{y\alpha}^{j}}{\widetilde
{F}_{y}^{j}}\right)  \right)  \\
&  =\scriptstyle\pm%
%TCIMACRO{\dsum \limits_{j=0}^{p-1}}%
%BeginExpansion
{\displaystyle\sum\limits_{j=0}^{p-1}}
%EndExpansion
\left(  \frac{\widetilde{F}_{x^{2}}^{j}-2\widetilde{F}_{xy}^{j}\widetilde{\nu
}_{j}+\widetilde{F}_{y^{2}}^{j}\widetilde{\nu}_{j}^{2}}{\widetilde{F}_{x}^{j}%
}\partial_{\alpha}f^{j}+\left(  \frac{\widetilde{F}_{y^{2}}^{j}}{\widetilde
{F}_{y}^{j}}-\frac{\widetilde{F}_{xy}^{j}}{\widetilde{F}_{x}^{j}}\right)
\frac{\widetilde{F}_{\alpha}^{j}}{\widetilde{F}_{y}^{j}}+\frac{\widetilde
{F}_{x\alpha}^{j}}{\widetilde{F}_{x}^{j}}-\frac{\widetilde{F}_{y\alpha}^{j}%
}{\widetilde{F}_{y}^{j}}\right)
\end{align*}
as desired.
\end{pf}

\begin{prop}
\label{third}The third derivative of $f^{k}$ defined using the system
(\ref{implic}) along the orbit is given by%
\begin{multline*}
\partial_{x^{3}}f^{k}=\scriptstyle\frac{\left(  \partial_{x^{2}}f^{k}\right)
^{2}}{\partial_{x}f^{k}}+\partial_{x}f^{k}%
%TCIMACRO{\dsum \limits_{j=0}^{k-1}}%
%BeginExpansion
{\displaystyle\sum\limits_{j=0}^{k-1}}
%EndExpansion
\frac{F_{x^{2}}^{j}-2F_{xy}^{j}\nu_{j}+F_{y^{2}}^{j}\nu_{j}^{2}}{F_{x}^{j}%
}\partial_{x^{2}}f^{j}\\
\scriptstyle+\partial_{x}f^{k}%
%TCIMACRO{\dsum \limits_{j=0}^{k-1}}%
%BeginExpansion
{\displaystyle\sum\limits_{j=0}^{k-1}}
%EndExpansion
\frac{F_{x^{3}}^{j}-3F_{x^{2}y}^{j}\nu_{j}+3F_{xy^{2}}^{j}\nu_{j}^{2}%
-F_{y^{3}}^{j}\nu_{j}^{3}}{F_{x}^{j}}\left(  \partial_{x}f^{j}\right)  ^{2}\\
\scriptstyle+\partial_{x}f^{k}%
%TCIMACRO{\dsum \limits_{j=0}^{k-1}}%
%BeginExpansion
{\displaystyle\sum\limits_{j=0}^{k-1}}
%EndExpansion
\left(  -\left(  F_{x^{2}}^{j}\right)  ^{2}+F_{x^{2}}^{j}F_{xy}^{j}\nu
_{j}+\left(  F_{x^{2}}^{j}{}F_{y^{2}}^{j}+2\left(  F_{xy}^{j}\right)
^{2}\right)  \nu_{j}^{2}-5F_{xy}^{j}F_{y^{2}}^{j}\nu_{j}^{3}+2\left(
F_{y^{2}}^{j}\right)  ^{2}\nu_{j}^{4}\right)  \left(  \frac{\partial_{x}f^{j}%
}{F_{x}^{j}}\right)  ^{2}%
\end{multline*}
with $\partial_{x}f^{j}$, $\partial_{x^{2}}f^{j}$ known from the previous results.

At the bifurcation point we obtain%
\begin{multline*}
\partial_{x^{3}}f^{p}=\scriptstyle\pm\left(  \partial_{x^{2}}f^{p}\right)
^{2}\pm%
%TCIMACRO{\dsum \limits_{j=0}^{p-1}}%
%BeginExpansion
{\displaystyle\sum\limits_{j=0}^{p-1}}
%EndExpansion
\frac{\widetilde{F}_{x^{2}}^{j}-2\widetilde{F}_{xy}^{j}\widetilde{\nu}%
_{j}+\widetilde{F}_{y^{2}}^{j}\widetilde{\nu}_{j}^{2}}{\widetilde{F}_{x}^{j}%
}\partial_{x^{2}}f^{j}\\
\scriptstyle\pm%
%TCIMACRO{\dsum \limits_{j=0}^{p-1}}%
%BeginExpansion
{\displaystyle\sum\limits_{j=0}^{p-1}}
%EndExpansion
\left(  \frac{\widetilde{F}_{x^{3}}^{j}-3\widetilde{F}_{x^{2}y}^{j}%
\widetilde{\nu}_{j}+3\widetilde{F}_{xy^{2}}^{j}\widetilde{\nu}_{j}%
^{2}-\widetilde{F}_{y^{3}}^{j}\widetilde{\nu}_{j}^{3}}{\widetilde{F}_{x}^{j}%
}\right)  \left(  \partial_{x}f^{j}\right)  ^{2}\\
\scriptstyle\pm%
%TCIMACRO{\dsum \limits_{j=0}^{p-1}}%
%BeginExpansion
{\displaystyle\sum\limits_{j=0}^{p-1}}
%EndExpansion
\left(  \frac{-\left(  \widetilde{F}_{x^{2}}^{j}\right)  ^{2}+\widetilde
{F}_{x^{2}}^{j}\widetilde{F}_{xy}^{j}\widetilde{\nu}_{j}+\left(  \widetilde
{F}_{x^{2}}^{j}\widetilde{F}_{y^{2}}^{j}+2\left(  \widetilde{F}_{xy}%
^{j}\right)  ^{2}\right)  \widetilde{\nu}_{j}^{2}}{\left(  \widetilde{F}%
_{x}^{j}\right)  ^{2}}+\frac{-5\widetilde{F}_{xy}^{j}\widetilde{F}_{y^{2}}%
^{j}\widetilde{\nu}_{j}^{3}+2\left(  \widetilde{F}_{y^{2}}^{j}\right)
^{2}\widetilde{\nu}_{j}^{4}}{\left(  \widetilde{F}_{x}^{j}\right)  ^{2}%
}\right)  \partial_{x}f^{j},
\end{multline*}

\end{prop}

\begin{pf}
We recall the second derivative from (\ref{SEC})
\begin{equation}
\partial_{x^{2}}f^{k}=\partial_{x}f^{k}%
%TCIMACRO{\dsum \limits_{j=0}^{k-1}}%
%BeginExpansion
{\displaystyle\sum\limits_{j=0}^{k-1}}
%EndExpansion
\frac{F_{x^{2}}^{j}\left(  F_{y}^{j}\right)  ^{2}-2F_{xy}^{j}F_{x}^{j}%
F_{y}^{j}+F_{y^{2}}^{j}\left(  F_{x}^{j}\right)  ^{2}}{F_{x}^{j}\left(
F_{y}^{j}\right)  ^{2}}\partial_{x}f^{j}. \label{second}%
\end{equation}
We have also%
\[
\partial_{x}f^{j+1}=-\frac{F_{x}^{j}}{F_{y}^{j}}\partial_{x}f^{j}=-\nu
_{j}\partial_{x}f^{j}\text{,}%
\]
and%
\begin{align*}
\partial_{x^{2}}f^{j+1}  &  =-\frac{F_{x}^{j}}{F_{y}^{j}}\partial_{x^{2}}%
f^{j}+\frac{F_{x^{2}}^{j}\left(  F_{y}^{j}\right)  ^{2}-2F_{xy}^{j}F_{x}%
^{j}F_{y}^{j}+F_{y^{2}}^{j}\left(  F_{x}^{j}\right)  ^{2}}{\left(  F_{y}%
^{j}\right)  ^{3}}\left(  \partial_{x}f^{j}\right)  ^{2}\\
&  =-\nu_{j}\partial_{x^{2}}f^{j}+\frac{F_{x^{2}}^{j}-2F_{xy}^{j}\nu
_{j}+F_{y^{2}}^{j}\nu_{j}^{2}}{F_{y}^{j}}\left(  \partial_{x}f^{j}\right)
^{2}\text{.}%
\end{align*}

The result is obtained differentiating (\ref{second}) and substituting
$\partial_{x}f^{j+1}$ and $\partial_{x^{2}}f^{j+1}$ by the expressions above
and simplifying. After some painful but straightforward computations we arrive
at the result.

At the bifurcation point we have $\partial_{x}f^{p}=\pm1$. Therefore, we get
easily the second statement.
\end{pf}

The classic Schwarzian derivative takes the form%
\[
Sf^{p}=\frac{\partial_{x}^{3}f^{p}}{\partial_{x}f^{p}}-\frac{3}{2}\left(
\frac{\partial_{x}^{2}f^{p}}{\partial_{x}f^{p}}\right)  ^{2}.
\]
In the case of implicitly defined dynamical systems, the Schwarzian derivative
can be computed using the previous results, giving%
\begin{multline*}
Sf^{k}\left(  x_{0}\right)  =\scriptstyle%
%TCIMACRO{\dsum \limits_{j=0}^{k-1}}%
%BeginExpansion
{\displaystyle\sum\limits_{j=0}^{k-1}}
%EndExpansion
\frac{F_{x^{2}}^{j}-2F_{xy}^{j}\nu_{j}+F_{y^{2}}^{j}\nu_{j}^{2}}{F_{x}^{j}%
}\partial_{x^{2}}f^{j}\\
\scriptstyle+%
%TCIMACRO{\dsum \limits_{j=0}^{k-1}}%
%BeginExpansion
{\displaystyle\sum\limits_{j=0}^{k-1}}
%EndExpansion
\frac{F_{x^{3}}^{j}-3F_{x^{2}y}^{j}\nu_{j}+3F_{xy^{2}}^{j}\nu_{j}^{2}%
-F_{y^{3}}^{j}\nu_{j}^{3}}{F_{x}^{j}}\left(  \partial_{x}f^{j}\right)  ^{2}\\
\scriptstyle+%
%TCIMACRO{\dsum \limits_{j=0}^{k-1}}%
%BeginExpansion
{\displaystyle\sum\limits_{j=0}^{k-1}}
%EndExpansion
\left(  -\left(  F_{x^{2}}^{j}\right)  ^{2}+F_{x^{2}}^{j}F_{xy}^{j}\nu
_{j}+\left(  F_{x^{2}}^{j}{}F_{y^{2}}^{j}+2\left(  F_{xy}^{j}\right)
^{2}\right)  \nu_{j}^{2}-5F_{xy}^{j}F_{y^{2}}^{j}\nu_{j}^{3}+2\left(
F_{y^{2}}^{j}\right)  ^{2}\nu_{j}^{4}\right)  \left(  \frac{\partial_{x}f^{j}%
}{F_{x}^{j}}\right)  ^{2}\\
\scriptstyle-\frac{1}{2}\left(
%TCIMACRO{\dsum \limits_{j=0}^{k-1}}%
%BeginExpansion
{\displaystyle\sum\limits_{j=0}^{k-1}}
%EndExpansion
\frac{F_{x^{2}}^{j}-2F_{xy}^{j}\nu_{j}+F_{y^{2}}^{j}\nu_{j}^{2}}{F_{x}^{j}%
}\partial_{x}f^{j}\right)  ^{2}\text{.}%
\end{multline*}
Although the rather long expression, the Schwarzian derivative can be easily
computed. In the case of the pitchfork, the last term vanishes.

Combining all the results in this section, we are able to study the
codimension one bifurcations of implicitly defined one-dimensional discrete
dynamical systems.

\section{Examples}

In this section we give examples for fold, transcritical, pitchfork and flip
bifurcations for periodic orbits of implicitly defined dynamical discrete
dynamical systems.%
%TCIMACRO{\FRAME{ftbpFU}{2.5797in}{2.5676in}{0pt}{\Qcb{Period three saddle
%orbit generated by the fold bifurcation when $f^{3}$ crosses the diagonal in
%three points.}}{\Qlb{FIG1}}{implicit1.eps}%
%width 2.5797in;  height 2.5676in;  depth 0pt;  original-width 3.6115in;
%original-height 3.5942in;  cropleft "0";  croptop "1";  cropright "1";
%cropbottom "0";  filename 'implicit1.eps';file-properties "XNPEU";}} }%
%BeginExpansion
\begin{figure}
[ptb]
\begin{center}
\includegraphics[
height=2.5676in,
width=2.5797in
]%
{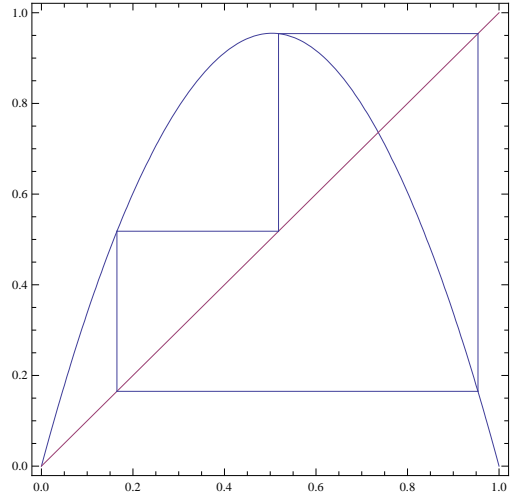}%
\caption{Period three saddle orbit generated by the fold bifurcation when
$f^{3}$ crosses the diagonal in three points.}%
\label{FIG1}%
\end{center}
\end{figure}
%EndExpansion
%TCIMACRO{\FRAME{ftbpFU}{2.5797in}{2.5676in}{0pt}{\Qcb{Triple tangency of the
%fold bifurcation for the implicit defined modified logistic.}}{\Qlb{FIG2}%
%}{implicit2.eps}{\special{ language "Scientific Word";  type "GRAPHIC";
%width 2.5797in;  height 2.5676in;  depth 0pt;  original-width 3.6115in;
%original-height 3.5942in;  cropleft "0";  croptop "1";  cropright "1";
%cropbottom "0";  filename 'implicit2.eps';file-properties "XNPEU";}} }%
%BeginExpansion
\begin{figure}
[ptb]
\begin{center}
\includegraphics[
height=2.5676in,
width=2.5797in
]%
{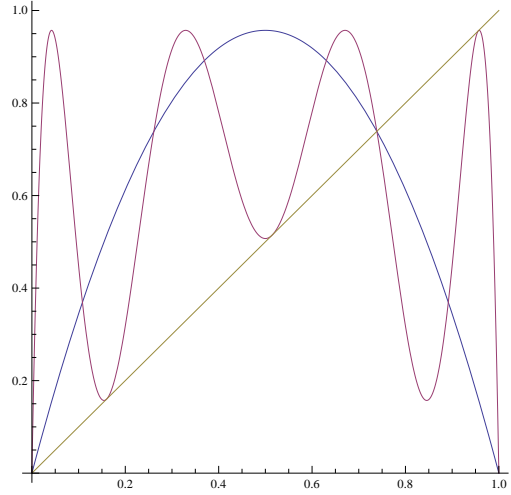}%
\caption{Triple tangency of the fold bifurcation for the implicit defined
modified logistic.}%
\label{FIG2}%
\end{center}
\end{figure}
%EndExpansion
%TCIMACRO{\FRAME{ftbpFU}{2.6187in}{2.5668in}{0pt}{\Qcb{The period two orbit at
%the transcritical bifurcation point. This non-hyperbolic orbit is a saddle,
%attracting from the outside and repelling to the inside of the interval.}%
%}{\Qlb{FIG5}}{implicit5.eps}{\special{ language "Scientific Word";
%type "GRAPHIC";  maintain-aspect-ratio TRUE;  display "USEDEF";
%valid_file "F";  width 2.6187in;  height 2.5668in;  depth 0pt;
%original-width 3.6115in;  original-height 3.5388in;  cropleft "0";
%filename 'implicit5.eps';file-properties "XNPEU";}} }%
%BeginExpansion
\begin{figure}
[ptb]
\begin{center}
\includegraphics[
height=2.5668in,
width=2.6187in
]%
{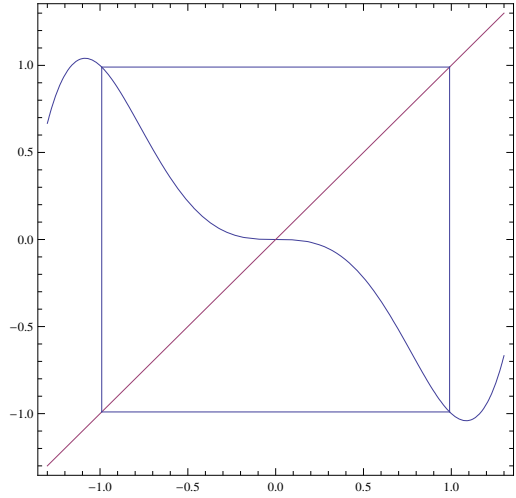}%
\caption{The period two orbit at the transcritical bifurcation point. This
non-hyperbolic orbit is a saddle, attracting from the outside and repelling to
the inside of the interval.}%
\label{FIG5}%
\end{center}
\end{figure}
%EndExpansion
%TCIMACRO{\FRAME{ftbpFU}{2.5668in}{2.5668in}{0pt}{\Qcb{Double tangency of the
%transcritical bifurcation of period $2$ for the implicit dynamical system of
%example \ref{ex2}. We can see $f^{2}$ for values of the parameter near $2$.
%The periodic points cross stabilities.}}{\Qlb{FIG6}}{implicit6.eps}%
%width 2.5668in;  height 2.5668in;  depth 0pt;  original-width 3.6115in;
%original-height 3.6115in;  cropleft "0";  croptop "1";  cropright "1";
%cropbottom "0";  filename 'implicit6.eps';file-properties "XNPEU";}} }%
%BeginExpansion
\begin{figure}
[ptb]
\begin{center}
\includegraphics[
height=2.5668in,
width=2.5668in
]%
{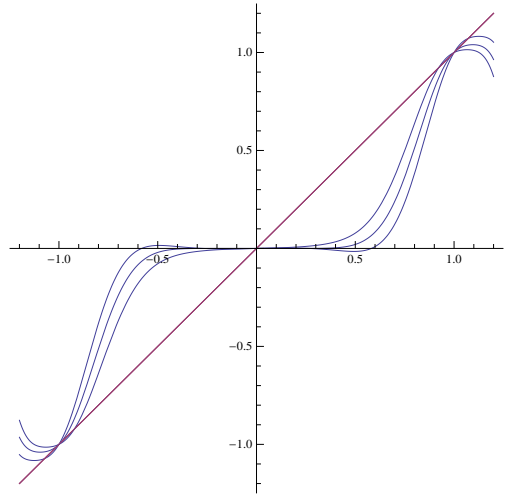}%
\caption{Double tangency of the transcritical bifurcation of period $2$ for
the implicit dynamical system of example \ref{ex2}. We can see $f^{2}$ for
values of the parameter near $2$. The periodic points cross stabilities.}%
\label{FIG6}%
\end{center}
\end{figure}
%EndExpansion
%TCIMACRO{\FRAME{ftbpFU}{2.6187in}{2.5668in}{0pt}{\Qcb{Non-hyperbolic orbit,
%although topologically stable, of period $2$ at the pitchfork bifurcation for
%the modified bimodal implicit dynamical system.}}{\Qlb{FIG3}}{implicit3.eps}%
%width 2.6187in;  height 2.5668in;  depth 0pt;  original-width 3.6115in;
%original-height 3.5388in;  cropleft "0";  croptop "1";  cropright "1";
%cropbottom "0";  filename 'implicit3.eps';file-properties "XNPEU";}} }%
%BeginExpansion
\begin{figure}
[ptb]
\begin{center}
\includegraphics[
height=2.5668in,
width=2.6187in
]%
{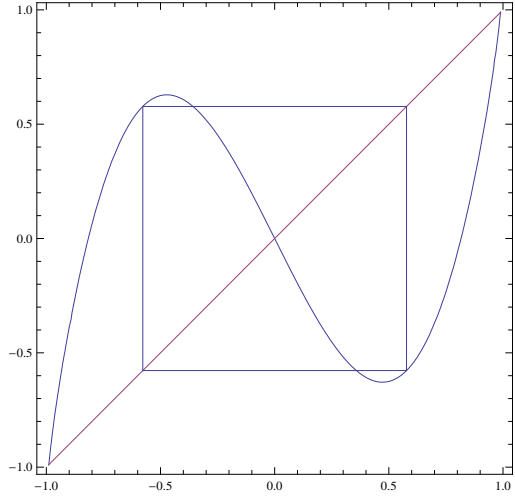}%
\caption{Non-hyperbolic orbit, although topologically stable, of period $2$ at
the pitchfork bifurcation for the modified bimodal implicit dynamical system.}%
\label{FIG3}%
\end{center}
\end{figure}
%EndExpansion
%TCIMACRO{\FRAME{ftbpFU}{2.6187in}{2.5668in}{0pt}{\Qcb{Double tangency of the
%pitchfork bifurcation of period $2$ for the implicit defined modified bimodal
%map. We can see $f$ and $f^{2}$. }}{\Qlb{FIG4}}{implicit4.eps}%
%width 2.6187in;  height 2.5668in;  depth 0pt;  original-width 3.6115in;
%original-height 3.5388in;  cropleft "0";  croptop "1";  cropright "1";
%cropbottom "0";  filename 'implicit4.eps';file-properties "XNPEU";}} }%
%BeginExpansion
\begin{figure}
[ptb]
\begin{center}
\includegraphics[
height=2.5668in,
width=2.6187in
]%
{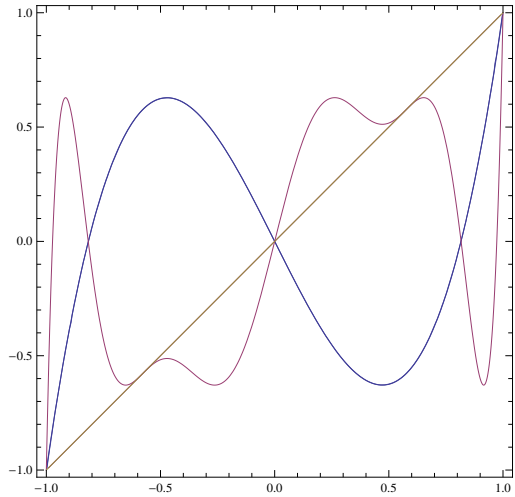}%
\caption{Double tangency of the pitchfork bifurcation of period $2$ for the
implicit defined modified bimodal map. We can see $f$ and $f^{2}$. }%
\label{FIG4}%
\end{center}
\end{figure}
%EndExpansion
%

%TCIMACRO{\FRAME{ftbpFU}{2.2628in}{2.2523in}{0pt}{\Qcb{Flip bifurcation. We see
%the double intersection of the map $f^{2}$ with the diagonal. The map has
%slope $-1$ at the relevant intersections.}}{\Qlb{FIG7}}{implicit7.eps}%
%width 2.2628in;  height 2.2523in;  depth 0pt;  original-width 3.3857in;
%original-height 3.3695in;  cropleft "0";  croptop "1";  cropright "1";
%cropbottom "0";  filename 'implicit7.eps';file-properties "XNPEU";}} }%
%BeginExpansion
\begin{figure}
[ptb]
\begin{center}
\includegraphics[
height=2.2523in,
width=2.2628in
]%
{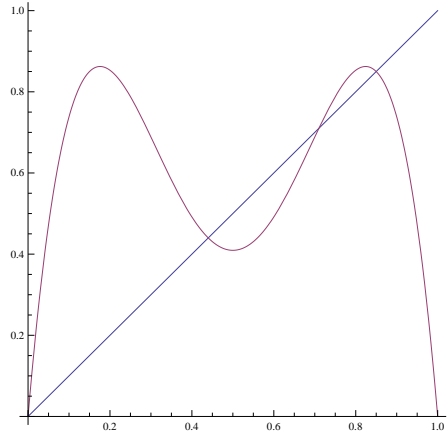}%
\caption{Flip bifurcation. We see the double intersection of the map $f^{2}$
with the diagonal. The map has slope $-1$ at the relevant intersections.}%
\label{FIG7}%
\end{center}
\end{figure}
%EndExpansion
%TCIMACRO{\FRAME{ftbpFU}{2.4647in}{2.4534in}{0pt}{\Qcb{Flip bifurcation. We can
%see the attracting period $4$ orbit generated by the flip bifurcation and in
%the center the repelling period $2$ orbit obtained from the original
%attracting period $2$ orbit.}}{\Qlb{FIG8}}{implicit8.eps}%
%width 2.4647in;  height 2.4534in;  depth 0pt;  original-width 3.3857in;
%original-height 3.3695in;  cropleft "0";  croptop "1";  cropright "1";
%cropbottom "0";  filename 'implicit8.eps';file-properties "XNPEU";}} }%
%BeginExpansion
\begin{figure}
[ptb]
\begin{center}
\includegraphics[
height=2.4534in,
width=2.4647in
]%
{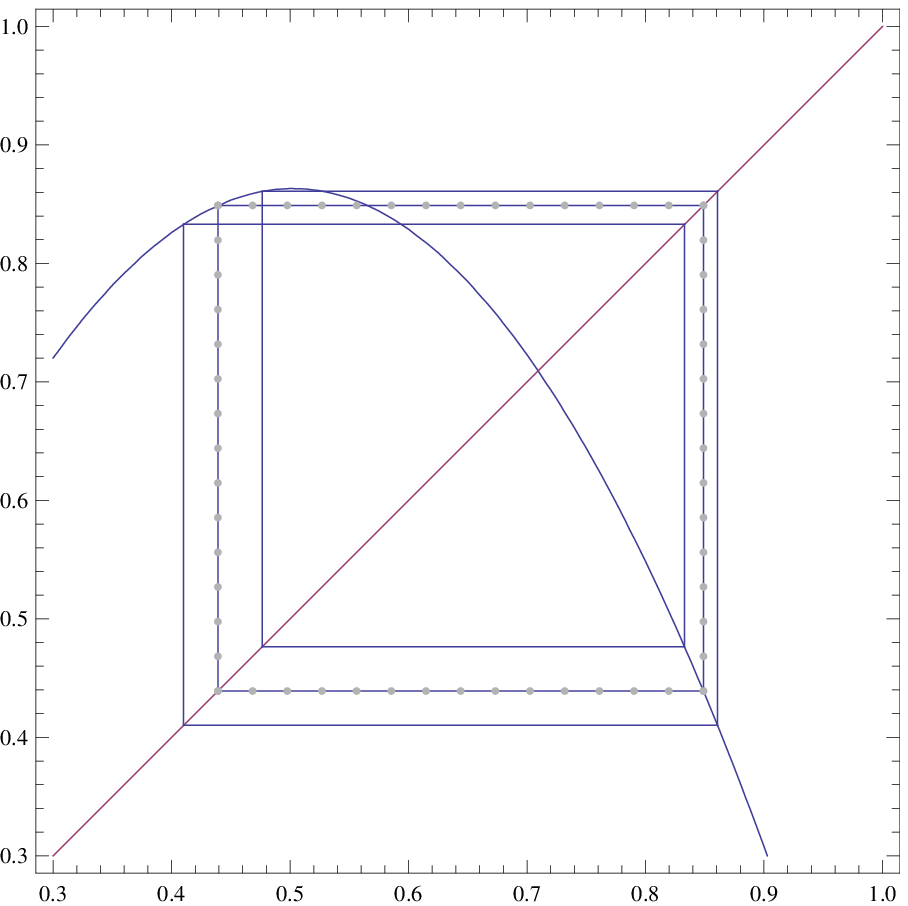}%
\caption{Flip bifurcation. We can see the attracting period $4$ orbit
generated by the flip bifurcation and in the center the repelling period $2$
orbit obtained from the original attracting period $2$ orbit.}%
\label{FIG8}%
\end{center}
\end{figure}
%EndExpansion

\begin{exmp}
\label{ex1}\textbf{Fold case, period 3. }Let be the implicitly defined
discrete dynamical system for $x_{n}\in\left[  0,1\right]  $ and $\alpha
\in\left[  0,4\right]  $, we call to the following model a modified implicit
logistic map%
\[
F\left(  x_{n},x_{n+1},\alpha\right)  =x_{n+1}-\alpha x_{n}(1-x_{n}%
+\frac{x_{n+1}^{P}}{B})=0.
\]
With $P=5$ and $B=100$, an explicit solution for $x_{n+1}$ is not possible,
since the map
\[
F\left(  x,y,\alpha\right)  =y-\alpha x(1-x+\frac{y^{5}}{100})=0,
\]
does not admit a closed formula for the solution $y$. The derivatives of $F$
are%
\begin{align*}
F_{x}\left(  x,y,\alpha\right)   &  =\alpha(-1+2x-\frac{y^{5}}{100}),\\
F_{y}\left(  x,y,\alpha\right)   &  =1-\frac{1}{20}\alpha xy^{4},\\
F_{\alpha}\left(  x,y,\alpha\right)   &  =-x(1-x+\frac{y^{5}}{100}),\\
F_{x^{2}}\left(  x,y,\alpha\right)   &  =2\alpha,\\
F_{y^{2}}\left(  x,y,\alpha\right)   &  =-\frac{\alpha xy^{3}}{5},\\
F_{xy}\left(  x,y,\alpha\right)   &  =-\frac{\alpha y^{4}}{20}.
\end{align*}
We are looking for a period $3$ fold, the bifurcation equations are%
\[
\left\{
\begin{array}
[c]{l}%
F\left(  x_{0},x_{1},\alpha\right)  =0,\\
F\left(  x_{1},x_{2},\alpha\right)  =0,\\
F\left(  x_{2},x_{0},\alpha\right)  =0,\\
\partial_{x}f^{3}\left(  x_{0}\right)  =\left(  -1\right)  ^{3}%
%TCIMACRO{\dprod \limits_{j=0}^{2}}%
%BeginExpansion
{\displaystyle\prod\limits_{j=0}^{2}}
%EndExpansion
\frac{F_{x}\left(  x_{j},x_{j+1\left(  \operatorname{mod}3\right)  }\right)
}{F_{y}\left(  x_{j},x_{j+1\left(  \operatorname{mod}3\right)  }\right)  }=1.
\end{array}
\right.
\]
A solution found numerically is%
\begin{align*}
x_{0}  &  =0.16498\ldots,x_{1}=0.51813\ldots,\\
x_{2}  &  =0.954\ldots,\alpha=3.75938\ldots.
\end{align*}
The previous computations show that there exists locally the implicitly
defined discrete dynamical system, since the derivative $F_{y}\left(
x,y,\alpha\right)  $ does not vanish in the the interval $\left[  0,1\right]
$ containing the orbit. The non-degeneracy condition holds at the periodic
orbit where $f=y$ is the implicitly defined iteration function%
\[
\partial_{x^{2}}f^{3}\left(  x_{0}\right)  =23.5\ldots\text{.}%
\]
The transversality condition gives
\[
\partial_{\alpha}f^{3}\left(  x_{0}\right)  =-0.844\ldots\text{.}%
\]
Therefore, the bifurcation is a supercritical fold with period three,
generating one period three attracting and one period three repelling orbits.
The saddle orbit at the bifurcation point can be seen in Figure \ref{FIG1}.
The bifurcation is via a simultaneous triple tangency at the diagonal, and can
be seen in Figure \ref{FIG2}.
\end{exmp}

\begin{exmp}
\label{ex2}\textbf{Transcritical case, period 2. }Let be the implicitly
defined discrete dynamical system for $x_{n}\in\left[  0,1\right]  $ and
$\alpha\in\left[  0,4\right]  $%
\[
F\left(  x_{n},x_{n+1},\alpha\right)  =x_{n+1}+x_{n}+\alpha x_{n}\left(
\left(  x_{n}-\frac{x_{n+1}^{P}}{B}\right)  ^{2}-1\right)  -x_{n}\left(
\left(  x_{n}-\frac{x_{n+1}^{P}}{B}\right)  ^{4}-1\right)  =0.
\]
With $P=3$ and $B=100$, an explicit solution for $x_{n+1}$ is not possible,
since the equation
\[
F\left(  x,y,\alpha\right)  =y+x-\alpha x\left(  \left(  x-\frac{y^{3}}%
{100}\right)  ^{2}-1\right)  -x\left(  \left(  x-\frac{y^{3}}{100}\right)
^{4}-1\right)  =0,
\]
does not admit a closed solution for $y$. At this point we omit the long
computations needed and the list of derivatives, for sake of brevity. The
reader can confirm our conclusions easily.

A period two solution for the transcritical bifurcation is found numerically
to be%
\[
x_{0}=0.9903\ldots,x_{1}=-0.9903\ldots,\alpha=2.
\]
There exists locally the implicitly defined discrete dynamical system, since
the derivative $F_{y}\left(  x,y,\alpha\right)  $ does not vanish in the the
interval $\left[  0,1\right]  $ containing the orbit. The non-degeneracy
condition holds at the periodic orbit%
\[
\partial_{x^{2}}f^{2}\left(  x_{0}\right)  =-16.79\ldots\text{,}%
\]
the derivative relative to the parameter is
\[
\partial_{\alpha}f^{2}\left(  x_{0}\right)  =0\text{.}%
\]
Therefore, the transversality condition is now%
\[
\partial_{\alpha x}f^{2}\left(  x_{0}\right)  =4.07769\text{.}%
\]
The conditions indicate a classical transcritical bifurcation, similar to the
one that happens for the logistic map at the origin, but for a period two
orbit. See Figures \ref{FIG5} and \ref{FIG6} for a graphical perspective of
this type of bifurcation.
\end{exmp}

\begin{exmp}
\label{ex3}\textbf{Pitchfork case, period 2. }Let be the implicitly defined
discrete dynamical system for $x_{n}\in\left[  0,1\right]  $ and $\alpha
\in\left[  0,4\right]  $, we call to the following model a modified implicit
bimodal map%
\[
F\left(  x_{n},x_{n+1},\alpha\right)  =x_{n+1}-\alpha\left(  x_{n}%
+\frac{x_{n+1}^{P}}{B}\right)  ^{3}-\left(  1-\alpha\right)  \left(
x_{n}+\frac{x_{n+1}^{P}}{B}\right)  =0.
\]
With $P=5$ and $B=100$, an explicit solution for $x_{n+1}$ is not possible,
since the map
\[
F\left(  x,y,\alpha\right)  =y-\alpha\left(  x+\frac{y^{5}}{B}\right)
^{3}-\left(  1-\alpha\right)  \left(  x+\frac{y^{5}}{B}\right)  =0,
\]
does not admit a closed formula for the solution $y$. We are looking for a
period $2$ pitchfork. The bifurcation equations are%
\[
\left\{
\begin{array}
[c]{l}%
F\left(  x_{0},x_{1},\alpha\right)  =0,\\
F\left(  x_{1},x_{0},\alpha\right)  =0,\\
\partial_{x}y^{2}\left(  x_{0}\right)  =\left(  -1\right)  ^{2}%
%TCIMACRO{\dprod \limits_{j=0}^{1}}%
%BeginExpansion
{\displaystyle\prod\limits_{j=0}^{1}}
%EndExpansion
\frac{F_{x}\left(  x_{j},x_{j+1\left(  \operatorname{mod}3\right)  }\right)
}{F_{y}\left(  x_{j},x_{j+1\left(  \operatorname{mod}3\right)  }\right)
}=1,\\
\partial_{x^{2}}y^{2}\left(  x_{0}\right)  =0.
\end{array}
\right.
\]
For sake of brevity we do not present here the derivatives of $F$ but only the
final results. A solution found numerically is%
\[
x_{0}=-0.5774599\ldots,x_{1}=0.5774599\ldots,\alpha=2.9989\ldots.
\]
meaning that there exists a periodic orbit with period two that bifurcates.
The first non-degeneracy condition holds at the periodic orbit%
\[
\partial_{x^{3}}y^{2}\left(  x_{0}\right)  =-295,6\ldots\text{.}%
\]
the first derivative in order to the parameter gives naturally%
\[
\partial_{\alpha}f^{2}\left(  x_{0}\right)  =0
\]
and the transversality condition is now%
\[
\partial_{\alpha x}f^{2}\left(  x_{0}\right)  =4.05\ldots\text{.}%
\]

Therefore, the bifurcation is a supercritical pitchfork with period two,
generating two new period two attracting orbits and the original period two
attracting orbit becomes repelling. The orbit at the bifurcation point can be
seen in Figure \ref{FIG3}. The bifurcation is via a simultaneous double
unfolding at the diagonal, and can be seen in Figure \ref{FIG4}.
\end{exmp}

\begin{exmp}
\label{ex4}\textbf{Flip case, period 2 into period 4. }Let be the implicitly
defined discrete dynamical system of example \ref{ex1} for $x_{n}\in\left[
0,1\right]  $ and $\alpha\in\left[  0,4\right]  $
\[
F\left(  x,y,\alpha\right)  =y-\alpha x(1-x+\frac{y^{5}}{100})=0,
\]
The new derivatives of $F$ that matter are%
\begin{align*}
F_{\alpha x}\left(  x,y,\alpha\right)   &  =-1+2x-\frac{y^{5}}{100},\\
F_{\alpha y}\left(  x,y,\alpha\right)   &  =-\frac{xy^{4}}{20},\\
F_{x^{3}}\left(  x,y,\alpha\right)   &  =0,\\
F_{x^{2}y}\left(  x,y,\alpha\right)   &  =0,\\
F_{xy^{2}}\left(  x,y,\alpha\right)   &  =-\frac{\alpha y^{3}}{5}\\
F_{y^{3}}\left(  x,y,\alpha\right)   &  =-\frac{3\alpha xy^{2}}{5}.
\end{align*}
We are looking for a period $2$ flip that bifurcates in a period $4$, the
bifurcation equations are for $x_{0}\not =x_{1}$%
\[
\left\{
\begin{array}
[c]{l}%
F\left(  x_{0},x_{1},\alpha\right)  =0,\\
F\left(  x_{1},x_{0},\alpha\right)  =0,\\
\partial_{x}f^{2}\left(  x_{0}\right)  =\left(  -1\right)  ^{2}%
%TCIMACRO{\dprod \limits_{j=0}^{1}}%
%BeginExpansion
{\displaystyle\prod\limits_{j=0}^{1}}
%EndExpansion
\frac{F_{x}\left(  x_{j},x_{j+1\left(  \operatorname{mod}3\right)  }\right)
}{F_{y}\left(  x_{j},x_{j+1\left(  \operatorname{mod}3\right)  }\right)  }=-1.
\end{array}
\right.
\]
A solution found numerically is%
\[
x_{0}=0.8466\ldots,x_{1}=0.4427\ldots,\alpha=3.405\ldots.
\]
The previous computations show that there exists locally the implicitly
defined discrete dynamical system, since the derivative $F_{y}\left(
x,y,\alpha\right)  $ does not vanish in the the interval $\left[  0,1\right]
$ containing the orbit. The non-degeneracy condition (\ref{Schwarz}) holds at
the periodic orbit where $f$ is the implicitly defined iteration function%
\begin{equation}
\frac{1}{2}\left(  f_{x^{2}}^{2}\left(  x_{0}\right)  \right)  ^{2}+\frac
{1}{3}f_{x^{3}}^{2}\left(  x_{0}\right)  =1383.1\not =0,
\end{equation}
which is equivalent to say that the Schwarzian derivative of $f^{2}$ is not
zero at $x_{0}$. The transversality condition (\ref{fliptransverse}) is%
\begin{equation}
f_{x\alpha}^{2}(x_{0})=1.45122\not =0\text{.}%
\end{equation}
Therefore, the bifurcation is a supercritical flip from period two to period
four, generating one period four attracting and one period two repelling
orbits. The bifurcation is via a simultaneous double $-1$ derivative for
$f^{2}$ at the diagonal, and can be seen in Figure \ref{FIG7}, finally in
Figure \ref{FIG8} we can see the orbits after the bifurcation, the dotted line
is the period two repelling orbit.
\end{exmp}

\subsubsection{Bifurcations in backward Euler and trapezoid methods}

Consider the autonomous differential equation%
\begin{equation}
x^{\prime}\left(  t\right)  =G\left(  x\left(  t\right)  \right)  \text{,
}x\left(  0\right)  =x_{0}\text{.} \label{diff}%
\end{equation}
In the usual Euler method the integral is estimated at the leftmost point of
each interval giving%
\begin{equation}
x_{n+1}-x_{n}=hG\left(  x_{n}\right)  , \label{Euler}%
\end{equation}
where $h$ is a positive real number, possibly very small. The backward, or
implicit Euler method \cite{Hai2009,Hai1996}, where the integral is estimated
using the rightmost point of each interval $x_{n+1}$ gives the iterative
scheme%
\begin{equation}
x_{n+1}-x_{n}=hG\left(  x_{n+1}\right)  . \label{Eulerb}%
\end{equation}
Actually this is a very simple one-dimensional discrete dynamical system,
obviously it can depend on internal parameters in $G$, but we are interested
in considering $h$ as the bifurcation parameter.

The iterative scheme is given by%
\begin{equation}
F\left(  x_{n},x_{n+1}\right)  =x_{n+1}-x_{n}-hG\left(  x_{n+1}\right)  =0.
\label{FEulerb}%
\end{equation}
Our function $F$ is%
\begin{equation}
F\left(  x,y\right)  =y-x-hG\left(  y\right)  .
\end{equation}
The original Euler method is considered explicit since in that case
$y=x+hG\left(  x\right)  $ and $x_{n+1}=x_{n}+hG\left(  x_{n}\right)  $.

In the case of the trapezoid method \cite{Hai2009} (which is a second order
method) we have for the same differential equation the iterative scheme%
\begin{equation}
F\left(  x_{n},x_{n+1}\right)  =x_{n+1}-x_{n}-\frac{h}{2}\left(  G\left(
x_{n+1}\right)  +G\left(  x_{n}\right)  \right)  =0 \label{Ftrapezoid}%
\end{equation}
and the function $F$ is%
\begin{equation}
F\left(  x,y\right)  =y-x-\frac{h}{2}\left(  G\left(  y\right)  +G\left(
x\right)  \right)  \text{,} \label{trapez}%
\end{equation}
this method is intrinsically implicit, since there is no immediate solution of
$F\left(  x,y\right)  =0$ for $y$. We consider now the existence of periodic
orbits in the Euler iteration, the period is $p$, the simplest case is the
asymptotic stable fixed point, which indicates that the solution of the
original differential equation has a limit when $t$ goes to infinity for a set
of initial conditions. Obviously the non-hyperbolic condition (\ref{b})\ for
the backward Euler method simplifies%
\[
\partial_{x}f^{p}\left(  x_{0}\right)  =\frac{1}{%
%TCIMACRO{\dprod \limits_{j=0}^{p-1}}%
%BeginExpansion
{\displaystyle\prod\limits_{j=0}^{p-1}}
%EndExpansion
\left(  1-hG^{\prime}\left(  x_{j+1\left(  \operatorname{mod}p\right)
}\right)  \right)  }=\pm1,
\]
this gives the non hyperbolic conditions for the backward Euler method%
\begin{equation}%
%TCIMACRO{\dprod \limits_{j=0}^{p-1}}%
%BeginExpansion
{\displaystyle\prod\limits_{j=0}^{p-1}}
%EndExpansion
\left(  1-hG^{\prime}\left(  x_{j}\right)  \right)  =\pm1\text{.}
\label{bEuler-nh}%
\end{equation}
For the trapezoid method the non-hyperbolic condition (\ref{b}) is%
\[%
%TCIMACRO{\dprod \limits_{j=0}^{p-1}}%
%BeginExpansion
{\displaystyle\prod\limits_{j=0}^{p-1}}
%EndExpansion
\frac{1+\frac{h}{2}G^{\prime}\left(  x_{j}\right)  }{1-\frac{h}{2}G^{\prime
}\left(  x_{j}\right)  }=\pm1\text{.}%
\]
The non-degeneracy condition (\ref{Dalphabifpoint}) is%

\[
\partial_{h}f^{p}\left(  x_{0}\right)  =\left(  -1\right)  ^{p}%
%TCIMACRO{\dsum \limits_{j=0}^{p-1}}%
%BeginExpansion
{\displaystyle\sum\limits_{j=0}^{p-1}}
%EndExpansion
\frac{\left(  -1\right)  ^{j}F_{h}\left(  x_{j},x_{j+1}\right)  }{F_{y}\left(
x_{j},x_{j+1}\right)  }%
%TCIMACRO{\dprod \limits_{i>j}^{p-1}}%
%BeginExpansion
{\displaystyle\prod\limits_{i>j}^{p-1}}
%EndExpansion
\frac{F_{x}\left(  x_{i},x_{i+1}\right)  }{F_{y}\left(  x_{i},x_{i+1}\right)
}\not =0.
\]
For the backward Euler method it gives%
\[%
%TCIMACRO{\dsum \limits_{j=1}^{p}}%
%BeginExpansion
{\displaystyle\sum\limits_{j=1}^{p}}
%EndExpansion%
%TCIMACRO{\dprod \limits_{i\geq j}^{p}}%
%BeginExpansion
{\displaystyle\prod\limits_{i\geq j}^{p}}
%EndExpansion
\frac{G\left(  x_{j}\right)  }{\left(  1-hG^{\prime}\left(  x_{i}\right)
\right)  }\not =0\text{, with }x_{p}=x_{0}.
\]
For the trapezoid method we have%
\[
\frac{1}{2}%
%TCIMACRO{\dsum \limits_{j=0}^{p-1}}%
%BeginExpansion
{\displaystyle\sum\limits_{j=0}^{p-1}}
%EndExpansion
\frac{\left(  -1\right)  ^{j}\left(  G\left(  x_{j}\right)  +G\left(
x_{j+1}\right)  \right)  }{1-\frac{h}{2}G^{\prime}\left(  x_{j+1}\right)  }%
%TCIMACRO{\dprod \limits_{i>j}^{p-1}}%
%BeginExpansion
{\displaystyle\prod\limits_{i>j}^{p-1}}
%EndExpansion
\frac{1+\frac{h}{2}G^{\prime}\left(  x_{i}\right)  }{1-\frac{h}{2}G^{\prime
}\left(  x_{i+1}\right)  }\not =0\text{, with }x_{p}=x_{0}.
\]
We study a simple example for the backward Euler method. Similar examples can
be constructed for the trapezoid method.

\begin{exmp}
Consider the simple differential equation%
\begin{equation}
x^{\prime}=x^{5}-1\text{, }x_{0}=0\text{.} \label{differential}%
\end{equation}
This equation can be solved by quadratures but it is impossible to obtain an
explicit expression for the solution. Applying Euler backward method we get%
\[
x_{n+1}-x_{n}=hx_{n+1}^{5}-h\text{,}%
\]
i.e.,
\[
F\left(  x_{n},x_{n+1}\right)  =x_{n+1}-x_{n}-h\left(  x_{n+1}^{5}-1\right)
=0\text{.}%
\]
Naturally, $G\left(  y\right)  =y^{5}-1$. We have $G^{\prime}\left(  y\right)
=-5y^{4}$. Equation (\ref{bEuler-nh}) is%
\begin{equation}%
%TCIMACRO{\dprod \limits_{j=0}^{p-1}}%
%BeginExpansion
{\displaystyle\prod\limits_{j=0}^{p-1}}
%EndExpansion
\left(  1-5hx_{j}^{4}\right)  =\pm1\text{.} \label{BackEulerH}%
\end{equation}
We have to solve (\ref{implic}) together with (\ref{BackEulerH}), we start by
the fixed point%
\[
\left\{
\begin{array}
[c]{l}%
x_{0}^{5}-1=0,\\
\left(  1-5hx_{0}^{4}\right)  =\pm1\text{.}%
\end{array}
\right.
\]
Excluding the trivial case $x=1$, $h=0$, there are no solutions for the fold
case. The solution is fairly simple for the flip case
\[
x_{0}=1\text{, }h=0.4\text{,}%
\]
This means that when $h=0.4$ the fixed point $x_{0}=1$ duplicates. When $h$ is
greater than $0.4$ the fixed point becomes attracting and is generated a
period two repelling orbit. Below $0.4$ the fixed point $x_{0}=1$ is repelling.

Now let us consider a period two orbit, the bifurcation equations are now%
\[
\left\{
\begin{array}
[c]{l}%
x_{1}-x_{0}=hx_{1}^{5}-h\\
x_{0}-x_{1}=hx_{0}^{5}-h\\
\left(  1-5hx_{0}^{4}\right)  \left(  1-5hx_{1}^{4}\right)  =\pm1\text{.}%
\end{array}
\right.
\]
We get, among complex solutions not considered here, the non trivial
($h\not =0$) real solutions for the fold case%
\begin{align*}
x_{0}  &  =x_{1}=1\text{, }h=0.4\text{, degenerate and obtained previously }\\
x_{0}  &  =1.15767\text{, }x_{1}=-0.602341\text{, }h=1.63071\text{.}%
\end{align*}
obviously $x_{0}=-0.602341$ and $x_{1}=1.15767$ is also a solution.

For the flip case we get the period doubling point where a period two orbit
duplicates its period%
\begin{align*}
x_{0}  &  =1.12579,\text{ }x_{1}=0.718620,\text{ }h=0.503700,\\
x_{0}  &  =-0.580682,\text{ }x_{1}=1.15618,\text{ }h=1.62930.
\end{align*}
This means that the previously created at $h=0.4$ repelling period two
solution, bifurcates again when $h=0.503700$ to a period $4$ orbit.

Finally, among other period three solutions, there is a period three fold at a
low value of $h$
\[
x_{0}=0.784072,\text{ }x_{1}=0.16453,\text{ }x_{2}=1.22008,\text{
}h=0.619616.
\]
Due to the continuity of all the functions involved this implies the existence
of chaos for low values of the parameters, even in the case of the backward
Euler method of a very simple first order differential equation. It is a well
known fact that one-dimensional discrete dynamical systems are more complex
than one-dimensional continuous dynamical systems. Nevertheless, the existence
of chaos for small values of the parameter $h$ is still exciting.
\end{exmp}

The previous example suggests the existence of a plethora of phenomena
deserving further research in implicit numeric methods. By force, the more
general cases of implicit discrete dynamical systems, which are very scarce in
the literature, are a vast field of research totally open.

\textbf{Acknowledgement} Partially funded by FCT/Portugal through
UID/MAT/04459/ 2013 for CMAGDS.

%\bibliographystyle{elsart-num-sort}
%\bibliography{BibloH2}

\end{document}